\documentclass[10pt]{article}

\usepackage{latexsym}
\usepackage{amsfonts}
\usepackage{amsmath}
\usepackage{mathrsfs}  
\usepackage{dsfont}    
\usepackage{bbold}     
\usepackage[english]{babel}
\usepackage{caption}
\usepackage{epsfig}
\usepackage{float}
\usepackage{subfigure}

\newcommand{\zerarcounters}{\setcounter{equation}{0}}

\newcommand{\ZZZ}{\mathds{Z}}

\newcommand{\NNN}{\mathds{N}}
\newcommand{\QQQ}{\mathds{Q}}
\newcommand{\RRR}{\mathds{R}}

\newcommand{\one}{\mathds{1}}

\newcommand{\gotD}{{\mathfrak D}}

\newcommand{\gotK}{{\mathfrak K}}

\newcommand{\io}{\infty}
\newcommand{\eps}{\varepsilon}

\newcommand{\n}{\nu}

\newcommand{\om}{\omega}

\newcommand{\s}{\sigma}

\newcommand{\wh}{\widehat}

\newcommand{\der}{{\rm d}}

\def\deriv#1#2#3{
\ifnum\catcode`#3=12    
   \ifnum#3=1
      \frac {{\rm d}#1} {{\rm d}#2}
   \else
      \frac {{\rm d}^{#3}#1} {{\rm d}#2^{#3}}
   \fi
\else   
      \frac {{\rm d}^{#3}#1} {{\rm d}#2^{#3}}
\fi}

\def\qed{\hfill\raise1pt\hbox{\vrule height5pt width5pt depth0pt}}

\def\ins#1#2#3{\vbox to0pt{\kern-#2 \hbox{\kern#1 #3}\vss}\nointerlineskip}

\newcommand{\Arnold}{Arnol$'\!$d}

\begin{document}

\title{\Arnold\ tongues for a resonant injection-locked frequency divider:
analytical and numerical results}
\author{
\bf M. V. Bartuccelli$^1$,
J. H.B. Deane$^1$,
G. Gentile$^2$,
F. Schilder$^1$
\vspace{2mm}
\\ \small 
$^1$Department of Mathematics,
University of Surrey, Guildford, GU2 7XH, UK.
\\ \small 
E-mails:
m.bartuccelli@surrey.ac.uk, j.deane@surrey.ac.uk, f.schilder@surrey.ac.uk
\\  \small  
$^2$Dipartimento di Matematica, Universit\`a di Roma Tre, Roma,
00146, Italy.
\\ \small
E-mail: gentile@mat.uniroma3.it 
}

\date{}

\maketitle

\begin{abstract}
In this paper we consider a resonant injection-locked frequency divider
which is of interest in electronics, and we investigate the
frequency locking phenomenon when varying the amplitude and frequency
of the injected signal. We study both analytically and
numerically the structure of the \Arnold\ tongues in the
frequency-amplitude plane. In particular, we provide exact
analytical formulae for the widths of the tongues, which correspond
to the plateaux of the devil's staircase picture.
The results account for numerical and experimental findings presented
in the literature for special driving terms and, additionally,
extend the analysis to a more general setting.
\end{abstract}




\zerarcounters
\section{Introduction}\label{sec:1}

The locking of oscillators onto subharmonics of the driving frequency
(known as \textit{frequency locking} or \textit{frequency
demultiplication}) has been well known in electronics since the work of
van der Pol and van der Mark~\cite{PM}; see also~\cite{KC}.
In the frequency-amplitude plane, locking occurs in distorted
wedge-shaped regions (\textit{\Arnold\ tongues}) with apices
corresponding to rational values on the (scaled) frequency axis.
If one plots the ratio of the driving frequency to
the output frequency versus the driving frequency,
one obtains a so-called \textit{devil's staircase}, i.e.
a self-similar fractal object, where the qualitative
structure is replicated at higher levels of resolution,
with plateaux corresponding to rational values of the ratio.

The frequency locking phenomenon, the existence of the \Arnold\ tongues,
and the devil's staircase structure have been proved rigorously in
some mathematical models, such as the circle map~\cite{A}, and studied
numerically or experimentally for several electronic circuits,
such as the van der Pol equation~\cite{GH,PL},
the Josephson junction~\cite{AC,L1,QWZ},
the Chua circuit~\cite{PZC} among others.

In this paper we are interested in studying both analytically and
numerically an electronic circuit, namely a resonant
\textit{injection-locked frequency divider} (ILFD), first considered
in~\cite{OBK}. In~\cite{OBYK} a differential equation was introduced
to describe the circuit and it was shown that the numerical
integration of the equation reliably reproduces the experimental data.

In~\cite{BDG3}, the differential equation describing the ILFD
was studied with the purpose of explaining analytically
the appearance of frequency locking. In particular,
the full differential equation in question was shown to be
of the form $u'' + u'h(u) + k(u) + \mu\Psi(u, u', t) = 0$,
where $h(u)$ and $k(u)$ are even and odd functions of $u$, respectively,
the prime denotes differentiation with respect to $t$, and $\mu$ is the
amplitude of the perturbation $\Psi$, which is taken to be periodic
in $t$, with frequency $\omega$, i.e. $\Psi(u, u', t)=
\Psi(u, u', t+2\pi/\omega)$. If the drive is purely sinusoidal,
as in~\cite{OBYK}, the Fourier series expansion of $\Psi(u, u', t)$,
$$ \Psi(u,u',t) = \sum_{\nu\in\ZZZ} {\rm e}^{{\rm i}\nu\omega t}
\Psi_{\nu}(u,u') , $$
contains only the first harmonics $\nu=\pm1$
(i.e. $\Psi_{\nu}(u,u')=0$ for $|\nu| >1$).

When $\mu = 0$, the system is unperturbed,
and the differential equation is of a particular form known as the
Li\'{e}nard equation~\cite{H,GH}; the best known example of this
type is the van der Pol equation. Under suitable assumptions on $h$ and
$k$, the Li\'{e}nard equation admits a globally-attracting limit cycle.

The phenomenon of frequency locking manifests itself in the ILFD when
the ratio of the driving frequency $\omega$
to the output frequency $\Omega$ is plotted against $\omega$.
When $\omega$ is close to a rational multiple $\rho=p/q$ of the
frequency $\Omega_{0}$ of the limit cycle of the unperturbed system,
then $\Omega$ is fixed such that $\omega=\rho\Omega$. Therefore
the  plot has a devil's staircase structure~\cite{OBK}, with plateaux
corresponding to rational values of the ratio $\omega/\Omega$.
If $\omega/\Omega=p/q$ one says that there is a \textit{resonance}
(or \textit{synchronisation}) of order $p\!:\!q$.
For purely sinusoidal perturbations $\Psi$, such as those considered
explicitly in~\cite{OBK,OBYK}, the main plateaux correspond to even
values of $\rho$ (a physical argument was given in~\cite{YXK}).
The perturbation theory approach taken in~\cite{BDG3}
successfully explains the experimental observations, by computing
quantitatively the way in which the widths of the plateaux depend
on the amplitude of the perturbation $\mu$, assumed small.

In an alternative visualisation of the situation, a two-dimensional
plot showing where locking takes place is constructed in the
$(\omega, \mu)$ plane. The \Arnold\ tongues have widths and centre-lines
which vary as some (explicitly computable) integer power
of $\mu$~\cite{BDG3}. The experimentally-observed
dominance of tongues for which the ratio $\omega/\Omega_{0}$ is close
to an even integer can be explained by the fact that only
these tongues have a width of order $\mu$: all other tongues grow
in width as some higher power of $\mu$. Specifically,
if $\rho\in\QQQ$ and $\Delta\omega(\rho)=\{ \om : \omega/\Omega =
\rho \}$ is the width of the corresponding \textit{locking interval}
at fixed $\mu$, it was proved that, for sinusoidal perturbations,
\begin{equation}
\Delta\omega(2n/k) = O(\mu^{k}) , \qquad
\Delta\omega((2n-1)/k) = O(\mu^{2k})
\label{eq:1.1} \end{equation}
for all $k,n \in \NNN$ such that $2n/k$ and $(2n +1)/k$ are
irreducible fractions. The centre-lines are vertical for $\rho = 2n$
and bend away from the vertical by a distance of order $\mu^2$
for all other values of $\rho$.  In~\cite{BDG3} we also stressed
that the property~(\ref{eq:1.1}) strongly depends on the
particular form of the drive, more precisely on the fact that,
as in~\cite{OBK,OBYK}, the driving was taken to contain only
the first harmonics.

More generally, one can express $\Delta\omega(\rho)$ as
a power series (\textit{perturbation series}) in $\mu$,
\begin{equation}
\Delta\omega(\rho) = \sum_{k=1}^{\io} \mu^{k} \Delta_{k}\omega(\rho) .
\label{eq:1.2} \end{equation}
If $k_{0}\in\NNN$ is the
first integer such that $\Delta\omega_{k_{0}}(\rho)\neq0$
then from (\ref{eq:1.2}) one obtains $\Delta\omega(\rho) =
\mu^{k_{0}} \Delta_{k_{0}}\omega(\rho) + O(\mu^{k_{0}+1})$.

The convergence of the perturbation series for $\mu$ small enough ---
yielding analyticity in $\mu$ in a neighbourhood of the origin ---
was discussed and proved in~\cite{BDG3}; see also~\cite{GBD3}. Hence,
by keeping only the lowest order terms means that in~(\ref{eq:1.1})
we are looking at the leading contributions, without making any
uncontrolled truncation. 
The coefficients $\Delta_{k}\omega(\rho)$ are
given in the form of suitable integrals. However, it is not possible
to reduce this computation to the integration of elementary functions,
because the integrands involve functions which are known only
numerically. Thus, the computation of the integrals requires some work,
which we also discuss in this paper. 

A first order analysis of the locking intervals (in the same spirit
as in~\cite{BDG3}) is also performed in~\cite{AMOQR},
where only sinusoidal perturbations are considered;
in particular no prediction is made for resonances $p\!:\!q$
with $p/q\notin2\NNN$, as this would require a higher order analysis.
More general perturbations are considered in~\cite{M},
where a different approach is followed. However, this involves
approximations which, ultimately, correspond to a first order analysis.
By contrast, the analysis performed in~\cite{GBD3} and further
developed here allows us to go to arbitrary perturbation orders,
with a control on the remainder. Thus, not only one can find an exact
analytical expression for the leading order of the locking interval
of any resonance, but in principle one can also compute any
locking interval within any desired accuracy. In~\cite{DDK},
the lack of accurate analytical methods to predict the locking range
was deplored: in our opinion our analysis, which makes
no approximation, fills this gap. Of course, for practical purposes,
the computation of the locking interval for any given resonance
requires solving numerically some integrals (which become increasingly
complicated as the perturbation order increases). It would be
desirable to have a formula for the locking interval in terms of the
parameters $\alpha$ and $\beta$ of the system, were one to exist;
we point out that in~\cite{AMOQR} an asymptotic formula is given
in the limit of $\alpha=\io$ and $\beta$ large.

In further detail, the motivation for the current paper, which completes
the analysis of~\cite{BDG3} and also concentrates on numerics
related to the ILFD problem, is as follows:
\begin{enumerate}
\item \label{item1}
To compute the coefficients of the powers of $\mu$
explicitly, at least for the lowest perturbation orders,
so as to give a quantitative expression for the width of the tongues,
for more general perturbations than those considered in~\cite{BDG3}.
\item \label{item2}
To investigate numerically how large $\mu$ can be for
the analysis, which is carried out under the assumption
that $\mu\ll 1$, to break down.
\item \label{atd1}
To compute numerically the \Arnold\ tongue diagram
in the $(\omega,\mu)$ plane in the case that the periodic part
of the perturbation contains only one frequency, $\omega$.
This allows us to obtain information for values of $\mu$
where perturbation theory does not apply. On the other hand,
for smaller values of $\mu$, the analytical results provide a check
on the reliability of the numerical analysis.
\item \label{atd2}
The same as~\ref{atd1}, but in the case that the perturbation
contains all integer multiples of $\omega$: it was argued
in~\cite{BDG3} that the width of all tongues would then
be proportional to $\mu$ and all the centre lines would be vertical.
In particular we want to determine the constant of proportionality,
i.e. the coefficient $\Delta_{1}\omega(\rho)$ in~(\ref{eq:1.2}),
and show that the higher the values of $p$ and $q$
in $\rho=p/q$, the lower the constant.
\end{enumerate}

The rest of the paper is organised as follows. In Section~\ref{sec:2},
we summarise definitions and lemmas from~\cite{BDG3} which are needed
in the remainder of the paper, and extend the analysis to more general
analytical driving, possibly containing all harmonics.
In Section~\ref{sec:4} we concentrate on analytical results
concerning the \Arnold\ tongues, by gathering together all information
which can be obtained to any order of perturbation theory.
In Section~\ref{sec:3}, we describe the algorithms used to
carry out the computations of the integrals appearing in the theory.
In Section~\ref{sec:5} we give and discuss the numerical results:
after checking that they agree with the theory where the latter applies
(small $\mu$), we investigate how large $\mu$ can be for the theoretical
predictions to be reliably used. Finally in Section~\ref{sec:6}
we draw some conclusions.

\zerarcounters
\section{Preliminary analytical results}
\label{sec:2}

We recall the results of~\cite{BDG3}, and extend them to more general
perturbations. Numbered lemmas which we refer to in this paper are
taken directly from~\cite{BDG3}, and all proofs are given there too.
Reference to~\cite{BDG3} is given only for proofs and technical details,
the discussion below being quite self-contained.

The system of ordinary differential equations that describes the ILFD
can be put into the form
\begin{equation}
u'' + u'\,h(u)  + k(u) + \mu \Psi(u,u',t) = 0 ,
\label{eq:2.1} \end{equation}
with
\begin{equation}
h(u) := 1 - \beta +  3 \beta u^{2} ,
\qquad k(u) := u \left( \alpha - \beta + \beta u^{2} \right) ,
\label{eq:2.2} \end{equation}
and
\begin{equation}
\Psi(u,u',t) :=  u' \left( 3 u^{2} - 1 \right) f(\omega t) +
u \left( u^{2} - 1 \right) \left( f( \omega t) +
\omega f'(\omega t) \right) ,
\label{eq:2.3} \end{equation}
where here and henceforth $f'$ denotes the derivative of $f$
with respect to its argument.
The case $f(t)=\sin t$ (and hence $f'(t)=\cos t$) was explicitly
considered in~\cite{BDG3}. More generally one can consider any
analytic $2\pi$-periodic function
\begin{equation}
f(t) = \sum_{\nu=1}^{\io} \hat f_{\nu} \sin \nu t , \qquad
| \hat f_{\nu} | \le \Phi\,{\rm e}^{-\xi |\nu|} ,
\label{eq:2.4} \end{equation}
where the bound on the Fourier coefficients $\hat f_{\nu}$
--- for suitable positive constants $\Phi$ and $\xi$ --- follows
from the analyticity assumption on $f$.
For simplicity we confine ourselves
to odd functions: considering functions whose Fourier expansion
contains also cosines would overwhelm the analysis
without shedding further light on the results.

For $\mu = 0$,~(\ref{eq:2.1}) reduces to
the \textit{Li\'enard equation}~\cite{C,H}
\begin{equation}
u'' + u' \, h(u) + k(u) = 0 , 
\label{eq:2.5} \end{equation}
which we refer to as the `unperturbed equation'.
In order for it to have a globally-attracting limit cycle encircling
the origin~\cite{H,Z} we require that $\alpha > \beta > 1$
(this corresponds to the region of the parameter plane
called \textit{design area} in~\cite{DDK}). In that case,
we designate $u_0(t)$ the solution to~(\ref{eq:2.5}) corresponding
to the limit cycle. Let $T_{0}$ be the period of $u_0(t)$ and
let $\Omega_{0}=2\pi/T_{0}$ be the corresponding frequency:
$\Omega_{0}$ depends solely on $\alpha$ and $\beta$.

The unperturbed equation is autonomous, hence it clearly has the
property that if $u_0(t)$ is a solution, then so is $u_0(t+T)$
for any constant $T$. Consequently, we can fix the origin of time so
that $u_0(0) = U_0 > 0$ and $u'_0(0) = 0$. This has the effect of
shifting the third argument of $\Psi$ by some time $t_{0}$,
so $\Psi(u, u', t)$ becomes $\Psi(u, u', t+t_0)$ in~(\ref{eq:2.1}).

We also note that the symmetry properties of $h(u)$ and $k(u)$
guarantee that $u_{0}(t)$ has the property
\begin{equation}
u_{0}(t + T_{0}/2) = -u_{0}(t) \qquad \forall t\in\RRR ,
\label{eq:2.6} \end{equation}
which in turn yields that the Fourier expansion of
$u_{0}(t)$ contains only odd harmonics (lemma 2.1).

It is convenient to rescale time
by defining $\tau = \omega t$ so that $\Psi$ now has period $2\pi$ in
its third argument. After rescaling, the differential equation becomes
\begin{equation} \ddot u + \frac{1}{\omega} \dot u\,h(u) +
\frac{1}{\omega^{2}} k(u)
+ \mu \bar \Psi(u,\dot u,\tau+\tau_{0}) = 0 ,
\label{eq:2.7} \end{equation}
where a dot denotes differentiation with respect to $\tau$,
$\tau_{0}=\omega t_{0}$, and we have defined
\begin{equation}
\bar \Psi(u,\dot u,\tau) = 
\frac{1}{\omega^2}\left[
\omega\dot u \left( 3 u^{2} - 1 \right) f(\tau) +
u\left( u^{2} - 1 \right)\left( f(\tau) +
\omega f'(\tau) \right) \right] .
\label{eq:2.8} \end{equation}
We have shown in~\cite{BDG3} that if $\omega$ is `close' to
$p \Omega_{0}/q$, where $p,q\in\NNN$ are relatively prime,
then the frequency $\Omega$ of the solution exactly equals $q\omega/p$:
the system is said to be locked into the $p\!:\!q$ resonance.
How close $\omega$ has to be to $p \Omega_{0}/q$ depends
on $\mu$ and on the resonance itself --- quantitative
investigation of this `closeness' is the aim of the present paper.

Let $\rho=p/q\in\QQQ$. For $\omega$ close to $\rho\Omega_{0}$ put
\begin{equation}
\frac{1}{\omega} = \frac{1}{\rho\Omega_{0}} + \eps(\mu,\tau_{0}),
\qquad\mbox{where}\qquad
\eps(\mu,\tau_{0}) = \sum_{k=1}^{\io}
\mu^{k} \eps_{k}(\tau_{0}) .
\label{eq:2.9}
\end{equation}
Unlike~\cite{BDG3}, for the sake of convenience, here we
make explicit the dependence of $\eps$ on $\tau_{0}$.
The perturbation calculation is then carried out by substituting
the expression~(\ref{eq:2.9}) for $\omega$ in~(\ref{eq:2.7})
and expanding in powers of $\mu$. This results in
\begin{equation}
H(u,\dot u,\ddot u,\mu) := H_{0}(u,\dot u,\ddot u) +
\sum_{k=1}^{\io} \mu^{k} H_{k} (u,\dot u,\tau+\tau_{0}) = 0 , 
\label{eq:2.10}
\end{equation}
where
\begin{subequations}
\begin{align}
\hskip-.2truecm
H_{0}(u,\dot u,\ddot{u}) & =
\ddot{u} + \frac{\dot u\,h(u)}{\rho\Omega_0} +
\frac{k(u)}{\rho^2\Omega_0^2} ,
\label{eq:2.11a} \\
\hskip-.2truecm
H_{1}(u,\dot u,\tau) & =
\eps_{1}(\tau_{0}) \left( \dot{u}\, h(u) +
\frac{2\, k(u)}{\rho\Omega_{0}} \right) + 
\frac{\dot u \left(3 u^{2} - 1 \right)}{\rho\Omega_{0}} f(\tau)
+ u\left( u^{2} - 1 \right)
\left( \frac{f(\tau)}{\rho^{2}\Omega_{0}^{2}} 
+ \frac{f'(\tau)}{\rho\Omega_{0}} \right) ,
\label{eq:2.11b} \\
\hskip-.2truecm
H_{k}(u,\dot u,\tau) & =
\eps_{k}(\tau_{0}) \left( \dot{u}\, h(u) +
\frac{2\, k(u)}{\rho\Omega_{0}} \right) +
\sum_{k_{1}+k_{2}=k} \eps_{k_{1}}(\tau_{0})\,
\eps_{k_{2}}(\tau_{0}) \, k(u) \nonumber \\
& + \eps_{k-1}(\tau_{0}) \left[ \dot u \left(3 u^{2} - 1 \right)
f(\tau) + u\left( u^{2} - 1 \right) \left(
\frac{2f(\tau)}{\rho\Omega_{0}} + f'(\tau) \right) \right]
\label{eq:2.11c} \\
& + \sum_{k_{1}+k_{2}=k-1} \eps_{k_{1}}(\tau_{0})\,
\eps_{k_{2}}(\tau_{0}) \, u \left( u^{2} -1 \right) f(\tau) ,
\qquad \qquad k \ge 2 , \nonumber
\end{align}
\label{eq:2.11} \end{subequations}
\vskip-.3truecm\noindent
where the last line of~(\ref{eq:2.11c}) is missing for $k=2$.

In order to carry out the perturbation calculation to first order,
we first write the 
unperturbed system in the form
\begin{equation}
\dot{u} = v, \qquad \dot{v} = -\frac{v\,h(u)}{\rho\Omega_0}
-\frac{k(u)}{\rho^2 \Omega_0^2} \equiv G(u, v)
\label{eq:2.12} \end{equation}
which has a unique $2\pi\rho$-periodic solution
$\left(u_{0}(\tau),v_{0}(\tau)\right)$
such that $v_{0}(0)=0$. The Wronskian matrix of
equation~(\ref{eq:2.12}) is
\begin{equation}
W(\tau) = \left( \begin{matrix}
w_{11}(\tau) & w_{12}(\tau) \\
\dot{w}_{11}(\tau) & \dot{w}_{12}(\tau) \end{matrix} \right) 
\label{eq:2.13} \end{equation}
and satisfies
\begin{equation}
\begin{cases}
\dot W(\tau) = M(\tau) \, W(\tau) , & \\ W(0) = \one , \end{cases}
\qquad M(\tau) = \left( \begin{matrix}
0 & 1 \\ {\displaystyle
\frac{\partial}{\partial u} G(u_{0}(\tau), v_{0}(\tau)) }
& \displaystyle{
\frac{\partial}{\partial v} G(u_{0}(\tau), v_{0}(\tau)) }
\end{matrix} \right) .
\label{eq:2.14} \end{equation}
Lemma~4.1 then states that a solution to equation~(\ref{eq:2.14})
is obtained by setting
\begin{equation}
w_{12}(\tau) := c_{2} \dot u_{0}(\tau) , \qquad
w_{11}(\tau) := c_{1} \dot u_{0}(\tau) \int_{\bar \tau}^{\tau}
{\rm d} \tau' \frac{{\rm e}^{-F(\tau')}}{\dot u_{0}^{2}(\tau')} , 
\label{eq:2.15} \end{equation}
where $F(\tau)$ is given by
\begin{equation}
F(\tau) := \frac{1}{\rho\Omega_0}\int_{0}^{\tau} {\rm d}\tau' \, 
h(u_0(\tau')) ,
\label{eq:2.16}
\end{equation}
the constant $\bar \tau \in (0,\pi\rho)$ is chosen so that
$\dot w_{11}(0)=0$, while the constants $c_{1}$ and $c_{2}$ 
are such that $W(0) = \one$ --- it is shown in~\cite{BDG3} that this
choice can always be made.

By defining $r_{1}:=\ddot u_{0}(0)$, as in~\cite{BDG3}, and substituting
this into~(\ref{eq:2.12}), we find that 
\begin{equation}
r_1 = -\frac{U_{0} \left(\alpha - \beta + 
\beta U_0^2\right)}{\rho^2\Omega_0^2} .
\label{eq:2.17} \end{equation}
By Remark~4.2 in~\cite{BDG3}, we have, additionally,
that $c_1 = -r_1$ and $c_2 = 1/r_1$.

We further define $f_0$ by $\rho\Omega_0 f_0 = \langle h\rangle$,
the mean value of $h\left(u_0(\tau)\right)$, 
so that $f_0 = F(2\pi\rho)/(2\pi\rho)$, and we write
$F(\tau) = f_0\tau + \tilde{F}(\tau)$, where $\tilde{F}(\tau)$ is a
$2\pi\rho$-periodic function with zero mean. By lemma~1.2 one has
$f_{0}>0$ (cf. also~\cite{C}).

Lemma~4.4 then states that there exist two $2\pi\rho$-periodic 
functions $a(\tau)$ and $b(\tau)$ such that
\begin{equation}
w_{11}(\tau) = a(\tau) + {\rm e}^{-f_{0}\tau} \, b(\tau) , \qquad
w_{12}(\tau) = c \, a(\tau) ,
\label{eq:2.18} \end{equation}
for a suitable constant $c$. In order to develop perturbation theory
for a $2\pi p$-periodic solution, with $p\in\NNN$, which continues
the unperturbed solution when $\mu \neq 0$, one writes
\begin{equation}
u(\tau) = u_{0}(\tau) + \sum_{k=1}^{\io} \mu^{k} u_{k}(\tau) ,
\label{eq:2.19} \end{equation}
where $u_{0}(\tau)$ has period $2\pi\rho$ (and hence frequency $1/\rho$).
We have shown in~\cite{BDG3} that there exist
$2\pi p$-periodic functions $u_{k}(\tau)$ such that
the perturbation series~(\ref{eq:2.19}) converges for $\mu$ small enough.
The functions $u_{k}(\tau)$ are recursively defined
(see equation (7.2) of~\cite{BDG3}) as
\begin{equation}
u_{k}(\tau) = w_{11}(\tau) \bar u_{k} + w_{12}(\tau) \bar v_{k} +
\int_{0}^{\tau} \der\tau' \, {\rm e}^{F(\tau')} \left[
w_{12}(\tau) w_{11}(\tau') - w_{11}(\tau) w_{12}(\tau') \right]
\Psi_{k}(\tau) ,
\label{eq:2.20} \end{equation}
with
\begin{equation}
\Psi_{k}(\tau) := \left[ - \sum_{k'=1}^{k} \mu^{k'}
H_{k'}(u(\tau),\dot u(\tau),\tau+\tau_{0}) +
G_{2}(u(\tau),\dot u(\tau)) \right]_{k} ,
\label{eq:2.21} \end{equation}
where
\begin{eqnarray}
& & G_{2}(u,v) := G(u,v) - G(u_{0}(\tau),v_{0}(\tau)) \nonumber \\
& & \qquad \qquad - \left( u - u_{0}(\tau) \right)
\frac{\partial}{\partial u} G(u_{0}(\tau), v_{0}(\tau)) -
\left( v - v_{0}(\tau)) \right)
\frac{\partial}{\partial v} G(u_{0}(\tau), v_{0}(\tau))
\label{eq:2.22} \end{eqnarray}
and the notations $[\cdot]_{k}$ means that we expand
$u(\tau)$ and $\dot u(\tau)$ according to~(\ref{eq:2.19})
and, after taking the Taylor series of the functions $H_{k'}$,
$k'=1,\ldots,k$, and $G_{2}$, we keep the coefficients
of all contributions proportional to $\mu^{k}$.
In~(\ref{eq:2.20}), the initial conditions $\bar u_{k}$ must be suitably
fixed (again we refer to~\cite{BDG3} for details), whereas
$\bar v_{k}$ can be set equal to zero (cf. remark 5.1 of~\cite{BDG3}). 

Considering first order in $\mu$, we obtain the first order
compatibility condition  that has to be satisfied if $u_{1}(\tau)$
is to be periodic, i.e. $\langle {\rm e}^{\tilde{F}} b\,
\Psi_1\rangle =0$, where $\Psi_1(\tau) = - H_1(u_0(\tau), v_0(\tau),
\tau + \tau_0)$. Expanding $f(\tau)$ according to~(\ref{eq:2.4})
and using~(\ref{eq:2.11b}), this gives
\begin{equation}
\eps_{1}(\tau_{0})\, A + \sum_{\nu=1}^{\io} \hat f_{\nu}
\sum_{j=1}^{3} \left[ B_{j1\nu} \cos \nu\tau_{0} +
B_{j2\nu} \, \sin \nu\tau_{0} \right] = 0 ,
\label{eq:2.23} \end{equation}
where
\begin{eqnarray}
A \!\!\! & := & \!\!\! \frac{1}{2\pi\rho} \int_{0}^{2\pi\rho}
{\rm d} \tau \, {\rm e}^{\tilde F(\tau)} b(\tau)
\left[ \dot u_{0}(\tau) \, h(u_{0}(\tau)) +
\frac{2}{\rho\Omega_{0}} k ( u_{0}(\tau)) \right] ,
\label{eq:2.24} \end{eqnarray}
and
\begin{subequations}
\begin{align}
B_{i1\nu} & := \frac{1}{2\pi p} \int_{0}^{ 2\pi p}
{\rm d} \tau \, \frac{K_{i}(\tau)}{\rho^{2}\Omega_{0}^{2}}
\sin \nu\tau , \quad i=1,2, \qquad
B_{31\nu} := \nu \rho \Omega_{0} B_{22\nu} ,
\label{eq:2.25a} \\
B_{i2\nu} & := \frac{1}{2\pi p} \int_{0}^{ 2\pi p}
{\rm d} \tau \, \frac{K_{i}(\tau)}{\rho^{2}\Omega_{0}^{2}}
\cos \nu\tau , \quad i=1,2, \qquad
B_{32\nu} := - \nu \rho \Omega_{0} B_{21\nu} ,
\label{eq:2.25b}
\end{align}
\label{eq:2.25} \end{subequations}
\vskip-.3truecm\noindent
with
\begin{equation}
K_{1}(\tau) = {\rm e}^{\tilde F(\tau)} b(\tau)
\, \rho\Omega_{0} \, v_{0}(\tau)
\left( 3 u_{0}^{2}(\tau) - 1 \right) ,
\qquad K_{2}(\tau) = {\rm e}^{\tilde F(\tau)} b(\tau)
\, u_{0}(\tau) \left( u_{0}^{2}(\tau) - 1 \right) .
\label{eq:2.26} \end{equation}
By setting $D_{1\nu}=-\left( B_{11\nu} + B_{21\nu} + B_{31\nu} \right)$
and $D_{2\nu} = - \left( B_{12\nu} + B_{22\nu} + B_{32\nu} \right)$,
(\ref{eq:2.23}) then becomes
\begin{equation}
\eps_{1}(\tau_{0}) = \frac{1}{A} \sum_{\nu=1}^{\io} \hat f_{\nu} \left(
D_{1\nu} \cos \nu\tau_{0} + D_{2\nu} \sin \nu\tau_{0} \right)
:= \gotD_{1}(\tau_{0}) .
\label{eq:2.27} \end{equation}
By construction $\eps_{1}$ has zero mean, so that either
it is a non-constant function or it identically vanishes.
For purposes of comparison with~\cite{BDG3}, in the following we
shall shorten $D_{11}=D_{1}$ and $D_{21}=D_{2}$, and also
$B_{ij1}=B_{ij}$, which are the only relevant constants when $f$
contains only the first harmonics $\nu=1$ in~(\ref{eq:2.4}).

It is shown in Appendix~B of~\cite{BDG3} that
$A = -r_{1}\rho\Omega_{0}$; hence, from~(\ref{eq:2.17}), 
\begin{equation}
A = \frac{U_0\left(\alpha - \beta + \beta U_0^2\right)}{\rho\Omega_0} ,
\label{eq:2.28} \end{equation}
which provides an obvious means to check the numerics --- by
calculating $A$ from~(\ref{eq:2.24}) and comparing with~(\ref{eq:2.28}).

In~\cite{BDG3} it is also shown how to go to higher orders; to any
order $k\geq1$ one finds the compatibility condition
$\langle {\rm e}^{\tilde{F}} b\,\Psi_{k}\rangle =0$, where
the function $\Psi_{k}(\tau)$ is given by~(\ref{eq:2.21}).

The compatibility condition leads to all orders to
equations like~(\ref{eq:2.27}), which now read
\begin{equation}
\eps_{k}(\tau_{0}) = \gotD_{k}(\tau_{0}) ,
\qquad k \geq 1 ,
\label{eq:2.29} \end{equation}
for suitable functions $\gotD_{k}$ --- strictly speaking in~\cite{BDG3}
only the case $f(t)=\sin t$ is explicitly discussed, but one can
easily work out the general case of $f$ an arbitrary analytic function
by following the same strategy. Note that, with respect to~\cite{BDG3},
here we have included the factor $1/A$ in the definition of
$\gotD_{k}(\tau_{0})$.

The width of the plateau corresponding to a given $\rho$
(i.e. to a given resonance $p\!:\!q$ such that $\rho=p/q$)
can then be expressed as follows. First one defines
\begin{equation}
\gotD(\tau_{0},\mu) = \!\!
\sum_{k=1}^{\io} \eps^{k} \gotD_{k}(\tau_{0}) , \qquad
\eps_{\rm max}(\rho) := \!\!\!
\max_{0 \le \tau_{0} \le 2\pi} \gotD(\tau_{0},\mu) , \qquad
\eps_{\rm min}(\rho) := \!\!\!
\min_{0 \le \tau_{0} \le 2\pi} \gotD(\tau_{0},\mu) .
\label{eq:2.30} \end{equation}
Then by setting
\begin{equation}
\omega_{\rm min}(\rho) := \frac{\rho\Omega_{0}}{1 +
\rho\Omega_{0} \, \eps_{\rm max}(\rho)} , \qquad
\omega_{\rm max}(\rho) := \frac{\rho\Omega_{0}}{1 +
\rho\Omega_{0} \, \eps_{\rm min}(\rho)} ,
\label{eq:2.31} \end{equation}
the plateau corresponding to $\rho$ is given by
\begin{equation}
\Delta\omega(\rho) := \omega_{\rm max}(\rho)-\omega_{\rm min}(\rho) =
\frac{\rho^{2}\Omega_{0}^{2}
\left( \eps_{\rm max}(\rho) - \eps_{\rm min}(\rho) \right)}{
(1 + \rho\Omega_{0} \, \eps_{\rm min}(\rho))
(1 + \rho\Omega_{0} \, \eps_{\rm max}(\rho))} .
\label{eq:2.32} \end{equation}
In other words, for $\omega\in[\omega_{\rm min}(\rho),
\omega_{\rm max}(\rho)]$, one has locking $\omega=\rho\Omega$,
if $\Omega$ denotes the frequency of the output signal.
For each such value of $\omega$ the initial phase $\tau_{0}$
gets fixed to a value $\tau_{0}^{*}$ such that
$1/\omega=1/\rho\Omega_{0}+\eps(\mu,\tau_{0}^{*})$,
according to~(\ref{eq:2.9}). 

When the function $\eps_{1}(\tau_{0})$ in~(\ref{eq:2.27}) does not
vanish, then, if one further assumes that the second derivative
of $\gotD_{1}$ is non-zero at the stationary points (where
the maximum and minimum are attained), the first order approximation
is adequate. In other words, in such a case one can approximate
\begin{equation}
\eps_{\rm max}(\rho) = \mu
\max_{0 \le \tau_{0} \le 2\pi} \gotD_{1}(\tau_{0}) + O(\mu^{2}) ,
\qquad \eps_{\rm min}(\rho) = \mu
\min_{0 \le \tau_{0} \le 2\pi} \gotD_{1}(\tau_{0}) + O(\mu^{2}) ,
\label{eq:2.33} \end{equation}
and hence
\begin{subequations}
\begin{align}
\omega_{\rm min}(\rho) & =
\rho\Omega_{0} \left( 1 - \rho\Omega_{0}
\mu \max_{0 \le \tau_{0} \le 2\pi} \gotD_{1}(\tau_{0})
\right) + O(\mu^{2}) ,
\label{eq:2.34a} \\
\omega_{\rm max}(\rho) & =
\rho\Omega_{0} \left( 1 - \rho\Omega_{0}
\mu \min_{0 \le \tau_{0} \le 2\pi} \gotD_{1}(\tau_{0}) \right)
+ O(\mu^{2}) ,
\label{eq:2.34b}
\end{align}
\label{eq:2.34} \end{subequations}
\vskip-.3truecm\noindent
which gives a plateau of width
\begin{equation}
\Delta\omega(\rho) = \mu \Delta_{1}\omega(\rho) + O(\mu^{2}) , \qquad
\Delta_{1}\omega(\rho) := \rho^{2}\Omega_{0}^{2} \left(
\max_{0 \le \tau_{0} \le 2\pi} \gotD_{1}(\tau_{0}) -
\min_{0 \le \tau_{0} \le 2\pi} \gotD_{1}(\tau_{0}) \right)
\label{eq:2.35} \end{equation}
centred `around' the value $\omega_{\rm c}(\rho)=\rho\Omega_{0}$.
Since the function $\eps_{1}(\tau_{0})$ has zero mean, this means that
the corresponding \Arnold\ tongue in the $(\omega,\mu)$ plane
emanates from the point $\omega_{\rm c}(\rho)$ of the $\omega$-axis
as a cone with axis along the vertical and
angle $\theta(\rho)=\theta_{1}(\rho)+\theta_{2}(\rho)$
such that
\begin{equation}
\tan\theta_{1}(\rho)= - \rho^{2}\Omega_{0}^{2}
\min_{0 \le \tau_{0} \le 2\pi} \gotD_{1}(\tau_{0}) ,
\qquad \tan\theta_{2}(\rho)=\rho^{2}\Omega_{0}^{2}
\max_{0 \le \tau_{0} \le 2\pi} \gotD_{1}(\tau_{0}) .
\label{eq:2.36} \end{equation}
If $f$ contains only one harmonic, say $\hat f_{\nu}=0$ for
$|\nu|>1$ in~(\ref{eq:2.4}), then $\theta_{1}(\rho)=\theta_{2}(\rho)$,
and $\max\gotD_{1}(\tau_{0})$ $=$ $A^{-1}\sqrt{D_{1}^{2}+D_{2}^{2}}$.
Note that in such a case the second derivative of $\gotD_{1}$ equals
$\pm\hat f_{1}/A$ when the first derivative vanishes.

\zerarcounters
\section{\Arnold\ Tongues: analytical results}
\label{sec:4}

\subsection{First order contributions}
\label{sec:4.1}

Let us consider the expression in~(\ref{eq:2.32}) for the leading
contribution to the width of the plateau when the first order
contribution does not vanish. Then we neglect the high order terms
and approximate
\begin{equation}
\Delta\omega(\rho) \approx \mu \rho^{2}\Omega_{0}^{2} Q(\rho) ,
\qquad\mbox{where}\qquad
Q(\rho) = \max_{0 \le \tau_{0} \le 2\pi} \gotD_{1}(\tau_{0}) -
\min_{0 \le \tau_{0} \le 2\pi} \gotD_{1}(\tau_{0}) ,
\label{eq:4.1} \end{equation}
Note that, to obtain $\gotD_{1}(\tau_{0})$ from~(\ref{eq:2.27}),
one must keep only the summands such that $\hat f_{\nu} \neq 0$.

By writing $B_{ij\nu}$ according to~(\ref{eq:2.25}), one uses
that the Fourier expansions of the functions $K_{i}$
contain only even harmonics (cf. Section 6 in~\cite{BDG3}), i.e.
\begin{equation}
K_{i}(\tau) = \sum_{\substack{\nu'\in\ZZZ \\ \nu' \; {\rm even}}}
{\rm e}^{{\rm i}\nu' \tau/\rho} \hat K_{i\nu'} =
\sum_{\nu'\in \ZZZ} {\rm e}^{{\rm i}2\nu' \tau/\rho} \hat K_{i(2\nu')} .
\label{eq:4.2} \end{equation}
Furthermore, as~(\ref{eq:2.26}) shows, the functions
$K_{i}$ are analytic and hence the corresponding
Fourier coefficients $\hat K_{i\nu'}$ decay exponentially, i.e.
for $i=1,2 $ and for all $\nu'\in\ZZZ$ one has
$|\hat K_{i\nu'}| \le \Gamma {\rm e}^{-\xi_{1}|\nu'|}$
for suitable positive constants $\Gamma$ and $\xi_{1}$.

\begin{figure}[htbp]
\centering
\includegraphics*[angle=0,width=3.0in]{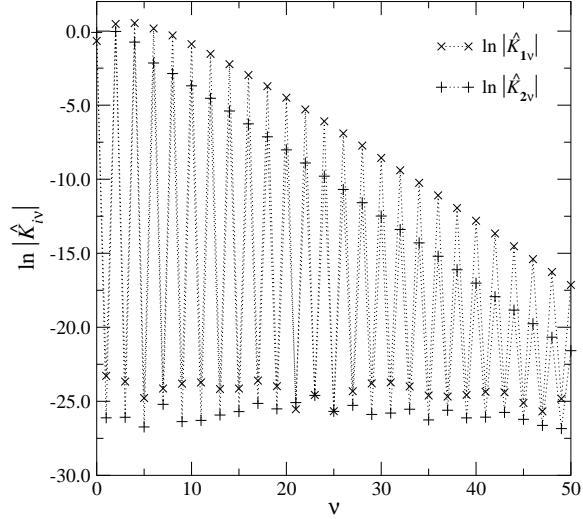}
\caption{Fourier coefficients of the functions $K_{1}(\tau)$, $\times$,
and $K_{2}(\tau)$, $+$, for $\alpha= 5$, $\beta = 4$. The odd
coefficients turn out to be zero (within the numerical error
of $\sim 10^{-11}$), according to (\ref{eq:4.2}), while all the
even coefficients are non-zero and decay exponentially.
The dotted lines are only to guide the eye.}
\label{fig:Kcoeffs}
\end{figure}

Hence by expanding $K_{i}$ according to~(\ref{eq:4.2}) and writing
\begin{equation}
\sin \nu\tau = \sum_{\sigma = \pm 1}
\frac{\sigma}{2{\rm i}} {\rm e}^{{\rm i}\sigma \nu\tau} , \qquad
\cos \nu\tau = \sum_{\sigma = \pm 1}
\frac{1}{2} {\rm e}^{{\rm i}\sigma \nu\tau} ,
\label{eq:4.3} \end{equation}
one realises that one can have $B_{ij\nu} \neq 0$ only if there exist
$\nu'\in\ZZZ$ such that $\hat K_{i\nu'} \neq 0$ and $\nu'q+\s\nu p = 0$.
If we assume that the first condition is satisfied for all
even $\nu'\in\ZZZ$ (numerical analysis ensures that such an assumption
is reasonable --- see figure~\ref{fig:Kcoeffs}), then
the key condition is that there exist $\nu'\in\ZZZ$ such that
\begin{equation}
2|\nu'|q=|\nu|p 
\label{eq:4.4} \end{equation}
with $p,q$ relatively prime integers. When this happens one has
\begin{equation}
B_{ij\nu} = \frac{1}{\rho^{2}\Omega_{0}^{2}}
\sum_{\substack{\nu'\in\ZZZ,\,\sigma=\pm 1 \\
2\nu'+\s\nu\rho =0}} \hat K_{i(2\nu')} \hat R_{j\s} ,
\quad i,j=1,2, \qquad\mbox{where}\qquad
\hat R_{1\s} = \frac{\s}{2{\rm i}} , \quad \hat R_{2\s} = \frac{1}{2} ,
\label{eq:4.5} \end{equation}
and $B_{31\nu}=\nu\rho\Omega_{0}B_{22\nu}$,
$B_{32\nu}=-\nu\rho\Omega_{0}B_{21\nu}$.
If the function $f$ contains only the first harmonics (so that
$\hat f_{\nu} \neq 0$ only for $|\n|=1$) then in~(\ref{eq:4.4})
one has to consider only the case $|\nu|=1$. Thus, as discussed
already in~\cite{BDG3}, one obtains $q=1$ and $p=2|\nu'|$,
i.e. $p$ must be even. This means that one finds plateaux of width
$O(\mu)$ only for resonances $p\!:\!q$ with $q=1$ and $p \in 2\NNN$.

On the contrary, if the function $f$ contains all the harmonics,
the condition~(\ref{eq:4.4}) has to be considered for all $\nu\in\ZZZ$,
and one finds easily non-vanishing contributions to~(\ref{eq:4.5}),
e.g. by taking $|\nu'|=p$ and $|\nu|=2q$. Thus, for any
resonance $p\!:\!q$ one has a plateau which to first order is
given by~(\ref{eq:4.1}). From~(\ref{eq:2.23}) one obtains
\begin{equation}
\gotD_{1}(\tau_{0}) = - \frac{1}{A \rho^{2}\Omega_{0}^{2} }
\sum_{\substack{\nu\geq 1 \\ \nu\rho \; {\rm even}}}
\hat f_{\nu} \left(
\overline K_{1(\nu\rho)} \,\cos \nu\tau_{0} +
\overline K_{2(\nu\rho)} \,\sin \nu\tau_{0} \right) ,
\label{eq:4.6} \end{equation}
where we have defined
\begin{subequations}
\begin{align}
\overline K_{1\nu} & := \sum_{\sigma=\pm 1}
\left[ \hat R_{1\sigma} 
\left( \hat K_{1(-\sigma\nu)} + \hat K_{2(-\sigma\nu)} \right) +
\nu \rho\Omega_{0} \hat R_{2\sigma}
\hat K_{2(-\sigma\nu)} \right] ,
\label{eq:4.7a} \\
\overline K_{2\nu} & := \sum_{\sigma=\pm 1} 
\left[ \hat R_{2\sigma}
\left( \hat K_{1(-\sigma\nu)} + \hat K_{2(-\sigma\nu)} \right) -
\nu \rho\Omega_{0} \hat R_{1\sigma}
\hat K_{2(-\sigma\nu)} \right] ,
\label{eq:4.7b}
\end{align}
\label{eq:4.7} \end{subequations}
\vskip-.3truecm\noindent
Let us define $\nu_{0}=\min\{ \nu \ge 1 : \nu\rho \; {\rm even}\}$
and $\nu_{1}=\min\{ \nu > \nu_{0} : \nu\rho \; {\rm even}\}$,
and set
\begin{equation}
\overline K_{\nu}(\rho) := \sqrt{ |\overline K_{1\nu}|^{2} +
|\overline K_{2\nu}|^{2}} , \qquad
Q_{0}(\rho) = \frac{2}{|A|\rho^{2}\Omega_{0}^{2}} | \hat f_{\nu_{0}} |
\, | \overline K_{\nu_{0}\rho}(\rho) | .
\label{eq:4.8} \end{equation}
Then one obtains
\begin{equation}
Q(\rho) = Q_{0}(\rho) + O
\left( \frac{|\hat f_{\nu_{1}}|\,|\overline K_{\nu_{1}\rho}(\rho)|}
{|\hat f_{\nu_{0}}|\,|\overline K_{\nu_{0}\rho}(\rho)|} \right) 
\label{eq:4.9} \end{equation}
which inserted into~(\ref{eq:4.1}) gives
\begin{equation}
\Delta\omega(\rho) \approx \frac{2\mu\rho\Omega_{0}}{|\bar r_{1}|}
| \hat f_{\nu_{0}} | \, | \overline K_{\nu_{0}\rho}(\rho) | ,
\label{eq:4.10} \end{equation}
where we have used that $A=-\bar r_{1}/\rho\Omega_{0}$,
with the constant $\bar r_{1}$ independent of $\rho\Omega_{0}$.  

If one keeps the whole sum in~(\ref{eq:4.6}) one finds,
always in the first order approximation,
\begin{equation}
\left| \Delta\omega(\rho) \right| \le
\frac{2 \mu \rho\Omega_{0}}{|\bar r_{1}|}
\max_{i=1,2} \sum_{\nu=1}^{\io} | \hat f_{\nu} | \,
| \overline K_{i(\nu\rho)}(\rho) | \le
\frac{2 \mu \rho\Omega_{0}}{|\bar r_{1}|} \max_{i=1,2}
\sum_{\substack{\nu\in\ZZZ \\ \nu\rho \; {\rm even}}}
\!\!\!\! \left( 1 + \nu \rho\Omega_{0} \right)
|\hat f_{\nu} | \, | \hat K_{i(\nu\rho)}(\rho) | .
\label{eq:4.11} \end{equation}
Since $p$ and $q$ are relatively prime the condition
$\nu\rho\in 2\ZZZ$ can be satisfied only if $|\nu|\geq q$
and $|\nu\rho| \geq p$. Therefore for fixed $\rho=p/q$ one has
\begin{equation}
\left| \Delta\omega(\rho) \right| \le \mu
C\, p^{2}q^{-1} {\rm e}^{-\xi_{1} p} {\rm e}^{-\xi q} ,
\label{eq:4.12} \end{equation}
where $C$ is a constant independent of $\rho$.
This shows that all the \Arnold\ tongues have width proportional
to $\mu$, but the constant of proportionality decays exponentially
with $p$ and $q$. Therefore, for fixed $\mu$, the union of all the
\Arnold\ tongues is $O(\mu)$, and hence tends to zero when $\mu\to0$,
as expected from~\cite{He}. 

For instance, if $f(t)=\sin t + \eta \sin 2 t$ in~(\ref{eq:2.4}),
one has $\Delta(2n)=c(2n)\,\mu + O(\mu^{2})$ and
$\Delta(2n-1)=c(2n-1)\eta\mu + O(\mu^{2})$,
for suitable constants $c(n)$ independent of $\mu$ and $\eta$.
Therefore, for all integer resonances the plateaux are of the same
order of magnitude --- provided, of course, $\eta$ is large
compared with $\mu$.

\subsection{Second order contributions}
\label{sec:4.2}

When the first order dominates, the second order gives a correction
which can be computed explicitly. When the first order vanishes,
the second order becomes the leading order (if it does not vanish too).

To compute the second order one needs the function $\gotD_{2}$
appearing in~(\ref{eq:2.29}) for $k=2$. The analysis in~\cite{BDG3}
and~(\ref{eq:2.29}) show that
\begin{equation}
\langle {\rm e}^{\tilde{F}} b\,\Psi_{2}\rangle =
A \eps_{2}(\tau_{0}) + \langle {\rm e}^{\tilde{F}} b\,\Xi_{2}
(\cdot;\tau_{0})\rangle \qquad \Longrightarrow \qquad
\gotD_{2}(\tau_{0}) = - \frac{1}{A} \langle {\rm e}^{\tilde F}
b \, \Xi_{2}(\cdot;\tau_{0}) \rangle ,
\label{eq:4.13} \end{equation}
where, by~(\ref{eq:2.11}) and~(\ref{eq:2.21}) with $k=2$,
one can write
\begin{subequations}
\begin{align}
\Xi_{2}(\tau;\tau_{0}) & = 
\widetilde\Xi_{2}(\tau;\tau_{0}) +
\overline\Xi_{2}(\tau;\tau_{0}) ,
\label{eq:4.14a} \\
\widetilde\Xi_{2}(\tau;\tau_{0}) & = -
\eps_{1}^{2}(\tau_{0}) \left[ \left(\alpha-\beta \right) u_{0}(\tau)
+ \beta u_{0}^{3}(\tau) \right] - \eps_{1}(\tau_{0}) v_{0}(\tau)
\left( 3 u_{0}^{2}(\tau)-1 \right) f(\tau+\tau_{0})
\nonumber \\
& - \eps_{1}(\tau_{0}) \left( u_{0}^{3}(\tau)-u_{0}(\tau) \right) \left(
\frac{2}{\rho\Omega_{0}} f(\tau+\tau_{0}) + f'(\tau+\tau_{0}) \right) ,
\label{eq:4.14b} \\
\overline\Xi_{2}(\tau;\tau_{0}) & = -
u_{1}(\tau) \frac{\partial H_{1}}{\partial
u_{0}}(u_{0}(\tau),v_{0}(\tau),\tau+\tau_{0}) -
\dot u_{1}(\tau) \frac{\partial H_{1}}{\partial
\dot u_{0}}(u_{0}(\tau),v_{0}(\tau),\tau+\tau_{0})
\nonumber \\
& + \frac{1}{2} u_{1}^{2}(\tau) \frac{\partial^{2} G}{\partial u^{2}}
(u_{0}(\tau),v_{0}(\tau)) +
u_{1}(\tau) \dot u_{1}(\tau) \frac{\partial^{2} G}{\partial u
\partial v}(u_{0}(\tau),v_{0}(\tau)) ,
\label{eq:4.14c} 
\end{align}
\label{eq:4.14} \end{subequations}
\vskip-.3truecm\noindent
with $H_{1}(u,v,\tau)$ and $G(u,v)$ given in~(\ref{eq:2.11b})
and~(\ref{eq:2.12}), respectively (we have explicitly used
that $G(u,v)$ is linear in $v$).

Thus, to compute~(\ref{eq:4.13}) one first needs the first order
solution $(u_{1},v_{1})$, with $v_{1}(\tau)=\dot u_{1}(\tau)$.
For $k=1$ equation ({\ref{eq:2.20})} gives
\begin{equation}
u_{1}(\tau) = c \, a(\tau)
\left( \mathcal{Q}_{2}(\tau) - \mathcal{Q}_{2}(0) -
\mathcal{Q}_{1}(0) \right) - c \, b(\tau)
\mathcal{Q}_{1}(\tau) ,
\label{eq:4.15} \end{equation}
where the functions $\mathcal{Q}_{1}$ and $\mathcal{Q}_{2}$
can  be read from equations~(5.3) and~(5.4) of~\cite{BDG3},
which we rewrite here for convenience,
\begin{equation}
\int_{0}^{\tau} {\rm d}\tau' {\rm e}^{F(\tau')}a(\tau')
\Psi_{1}(\tau') = {\rm e}^{f_{0}\tau}\mathcal{Q}_{1}(\tau) -
\mathcal{Q}_{1}(0) , \quad
\int_{0}^{\tau} {\rm d}\tau' {\rm e}^{\tilde F(\tau')}b(\tau')
\Psi_{1}(\tau') = \mathcal{Q}_{2}(\tau) -
\mathcal{Q}_{2}(0) ,
\label{eq:4.16} \end{equation}
and we are using that $\mathcal{Q}_{0}:=\langle {\rm e}^{\tilde F}
b\Psi_{1}\rangle=0$ and $\bar u_{1} = - c \mathcal{Q}_{1}(0)$.
Note that the functions $u_{1}$ and $v_{1}$
depend also on $\tau_{0}$; more precisely, by construction one has
\begin{equation}
u_{1}(\tau) = \sum_{\substack{\nu\in\ZZZ \\ \nu \; {\rm odd}}}
\sum_{\nu_{1}\in\ZZZ}
{\rm e}^{{\rm i}\nu \tau/\rho} {\rm e}^{{\rm i}\nu_{1}(\tau+\tau_{0})}
\hat U_{1\nu\nu_{1}} , \qquad
\hat U_{1\nu\nu_{1}} \propto \hat f_{\nu_{1}} ;
\label{eq:4.17} \end{equation}
as easily follows by reasoning as in the proof of lemma~8.2
in~\cite{BDG3}, the only difference being that $f$ can contain
all harmonics.

In particular, when $\eps_{1}$ vanishes
identically, then $\widetilde\Xi_{2}$ also is identically zero
and~(\ref{eq:4.13}) reduces to
\begin{eqnarray}
\gotD_{2}(\tau_{0}) 
& \!\!\! = \!\!\! &
\frac{1}{A} \frac{1}{2\pi p} \int_{0}^{2\pi p}
{\rm d}\tau \, {\rm e}^{\tilde F(\tau)} b(\tau)
\frac{1}{\rho\Omega_{0}} \Big\{ \Big[ \Big(
6 u_{0}(\tau) \, v_{0}(\tau)
+ \frac{1}{\rho\Omega_{0}}
\left( 3 u_{0}^{2}(\tau) - 1 \right) \Big)
f ( \tau + \tau_{0} ) \nonumber \\
& \!\!\! + \!\!\! &
\left( 3 u_{0}^{2}(\tau) - 1 \right)
f'( \tau + \tau_{0}) \Big] u_{1}(\tau) +
\left( 3 u_{0}^{2}(\tau) - 1 \right)
f( \tau + \tau_{0}) \, v_{1}(\tau)
\label{eq:4.18} \\
& \!\!\! + \!\!\! &
\frac{1}{2} \Big( v_{0}(\tau) h''(u_{0}(\tau)) +
\frac{k''(u_{0}(\tau))}{\rho\Omega_{0}} \Big) u_{1}^{2}(\tau)
+ h'(u_{0}(\tau)) \, u_{1}(\tau)\,v_{1}(\tau) \Big\} ,
\nonumber \end{eqnarray} 
where $h'(u)=6\beta u$, $h''(u)=6\beta$, and $k''(u)=6\beta u$ (here, as
well as for $f$, the prime denotes derivative with respect to the argument).

By using the expansion~(\ref{eq:4.17}) for $u_{1}$, one finds
that the function
\begin{equation}
\gotD_{2}(\tau_{0}) =
\sum_{\nu\in\ZZZ} {\rm e}^{{\rm i}\nu \tau_{0}} \gotD_{2,\nu} =
\gotD_{2,0} + \sum_{\substack{\nu\in\ZZZ \\ \nu \neq 0}}
{\rm e}^{{\rm i}\nu \tau_{0}} \gotD_{2,\nu} 
\label{eq:4.19} \end{equation}
can be written in the form
\begin{equation}
\gotD_{2}(\tau_{0}) = \frac{1}{2\pi p} \!\!
\sum_{\substack{\nu'\in\ZZZ \\ \nu' \; {\rm even}}} \!\!
\sum_{\nu_{1},\nu_{2}\in\ZZZ}
\int_{0}^{2\pi p} \!\!\!\!\!\!\!\! {\rm d}\tau 
\, {\rm e}^{{\rm i}\nu' \tau/\rho}
{\rm e}^{{\rm i}(\nu_{1}+\nu_{2})(\tau+\tau_{0})}
\hat\gotK_{\nu'\nu_{1}\nu_{2}} , \qquad
\hat\gotK_{\nu'\nu_{1}\nu_{2}} \propto
{\rm e}^{-\xi_{1}|\nu'|} \hat f_{\nu_{1}} \hat f_{\nu_{2}} ,
\label{eq:4.20} \end{equation}
for suitable coefficients $\hat\gotK_{\nu'\nu_{1}\nu_{2}}$. Then one 
sees that only the coefficients $\hat\gotK_{(2\nu')\nu_{1}\nu_{2}}$ with 
\begin{equation}
2|\nu'|q=|\nu_{1}+\nu_{2}|p , \qquad \hat f_{\nu_{1}}
\hat f_{\nu_{2}} \neq 0 ,
\label{eq:4.21} \end{equation}
contribute to~(\ref{eq:4.20}). The term with $\nu_{1}+\nu_{2}=\nu'=0$
gives the mean $\gotD_{2,0}$ of $\gotD_{2}$, and requires no condition
on $\rho$. This explains why the boundaries of the locking region
are either $O(\mu)$ --- when the first order dominates ---
or $O(\mu^{2})$ --- in all the other cases. However, the width of the
plateau arises from the variations of $\gotD_{2}$, hence it is related
to the terms in~(\ref{eq:4.20}) with $\nu\neq 0$ such that
(\ref{eq:4.21}) is satisfied. If there are any of such terms,
then the function $\gotD_{2}$ is not identically constant, and
therefore, in such a case, one has
\begin{equation}
\Delta\omega(\rho) = \mu^{2} \Delta_{2}\omega(\rho) + O(\mu^{3}) , \qquad
\Delta_{2}\omega(\rho) := \rho^{2}\Omega_{0}^{2} \left(
\max_{0 \le \tau_{0} \le 2\pi} \gotD_{2}(\tau_{0}) -
\min_{0 \le \tau_{0} \le 2\pi} \gotD_{2}(\tau_{0}) \right)
\label{eq:4.22} \end{equation}
which replaces~(\ref{eq:2.35}) when the first order vanishes.
For instance if $f$ contains only the first harmonics then
(\ref{eq:4.21}) is satisfied for $q=1$, $p\in\NNN$ and $\nu_{1}=
\nu_{2}=\pm1$, which shows that the plateaux corresponding to
odd $\rho$ are of order $\mu^{2}$ --- see~\cite{BDG3} for further details.

\subsection{Higher order contributions}
\label{sec:4.3}

If one wants to determine the higher order contributions,
the analysis above can be easily extended, even if it becomes
much more complicated from the computational point of view.
In general, if one writes
\begin{equation}
\Delta\omega(\rho) = \sum_{k=1}^{\io} \mu^{k} \Delta_{k}\omega(\rho) ,
\label{eq:4.23} \end{equation}
one finds
\begin{equation}
\Delta_{k}\omega(\rho) = \sum_{\nu\in\ZZZ} \!\!\!\!
\sum_{\substack{\nu_{1},\ldots,\nu_{k}\in\ZZZ \\
|\nu_{1}+\ldots+\nu_{k}|p=2\nu q}} \!\!\!\!\!\!\!\!\!
\Delta_{k}\omega(\rho;\nu_{1},\ldots,\nu_{k}) , \qquad
\Delta_{k}\omega(\rho;\nu_{1},\ldots,\nu_{k}) \propto
{\rm e}^{-2\xi_{1}|\nu|/\rho} \prod_{i=1}^{k} \hat f_{\nu_{i}} ,  
\label{eq:4.24} \end{equation}
so that, in order to single out the leading contribution to the width
of the plateau, one has to compare the size of the perturbation parameter
$\mu$ with the amplitudes of the harmonics $\hat f_{\nu}$ of the drive.

Note that to all orders $k$ the coefficients
$\Delta_{k}\omega(\rho)$ decay exponentially in both $p$ and $q$.
Thus, every time the first order does not vanish it
dominates, provided $\mu$ is small enough. If on the contrary
$\Delta_{k'}\omega(\rho) =0$ for all $1\le k' <k$ and
$\Delta_{k}\omega(\rho)\neq0$ then one has $\Delta\omega(\rho)=
O(\mu^{k})$ for $\mu$ small enough.

\zerarcounters
\section{Numerical computations}
\label{sec:3}

\subsection{Numerical solution of the ODE}
\label{nsol}

Since in general no closed-form solution to~(\ref{eq:2.7}) with $\mu = 0$
exists for $\beta\neq 0$, it is clear that this equation
must be solved numerically. Furthermore, in order to
approximate $u_0(\tau)$ and $\dot{u}_0(\tau)$, the ODE must be solved for a 
sufficiently long time that the transient has, for practical purposes, 
decayed to zero.  An effective initial procedure was found to be (i) 
solve the ODE from $\tau = -\tau_1$ to $\tau = 0$, where $\tau_1$ is
large, using any standard method,  for example, the Runge-Kutta
fourth order method; (ii) solve for a further small
time $\tau_2$ which is such that $\dot{u}_0(\tau_2) = 0$
and $u_0(\tau_2) > 0$, again using the Runge-Kutta method, 
and additionally using bisection to find $\tau_2$
such that the first condition is met; (iii) solve from $\tau = \tau_2$
to $\tau_3$, where $\tau_3$ is the smallest value of $\tau$ which is
greater than $\tau_2$ and for which, again, $\dot{u}_0(\tau_3) = 0$
and $u_0(\tau_3)>0$.
Then an estimate of the period of $u_0(\tau)$, $T_0$, is $\tau_3 -
\tau_2$ and an estimate of $U_0$ is $u_0(\tau_2)\approx u_0(\tau_3)$.

In practice, these estimates can then be somewhat improved by solving
the ODE assuming that a power series for $u_0(\tau)$ exists,
and computing this series around $\tau = 0$, using the initial 
conditions $u_0(0) = U_0$, as estimated
above, and $\dot{u}_0(0) = 0$. We can shift the origin of 
time from $\tau_3$ to zero because the ODE is autonomous.
Typically, several power series need to be
computed to cover the range $\tau = 0$ to $T_0$, but the
method has at least two advantages over Runge-Kutta. 
The first is that the error can be estimated by
implementing a test on the coefficients of the power series, as set
out in~\cite{Cha}; the second is that the Newton-Raphson method
can be used directly on the power series for the solution
around $T_0$ to find the value of $\tau$ for which
$\dot{u}_0(\tau) = 0$, and hence to estimate $T_0$.
The series used in practical computations had degree 30.

Once accurate values of $U_0$ and $T_0$ have been computed, it is a simple
matter to calculate a table of values of $u_0(\tau)$ and $\dot{u}_0(\tau)$ 
at $\tau = ih, i = 0\ldots M-1$ for some $M\in\mathds N$ and
for $h = T_0/M>0$ a given time-step. 

\subsection{Interpolation}
\label{interpol}

A discussion of a suitable interpolation method is now in order. 

In what follows, we will need to integrate functions of $u_0(\tau),
\dot{u}_0(\tau)$ and to do this we use an interpolation scheme from
which such integrals can be computed directly.

We start by discussing a scheme for interpolating from the values
of $u_0(\tau), \dot{u}_0(\tau)$ at discrete times $ih, i = 0\ldots
M-1$, produced by the numerical ODE solver.

Since $u_0(\tau)$ is periodic, the interpolation scheme
should reflect this --- standard methods based on the
Lagrange formula or Chebyshev polynomial interpolation 
are therefore not suitable. Instead, interpolation based on the function
\begin{equation}
I_K(\tau) = \frac{\sin K\pi \tau}{K\sin \pi \tau} ,
\label{I_of_t}
\end{equation}
where $K$ is an odd, positive integer, is used. This function is 
equivalent to a truncated Fourier series (see Appendix A) and has 
the properties that
(i) $I_K(\tau+1) = I_K(\tau)$, so it is periodic
(if $K$ is even, the period is not 1 but 2); and (ii) 
$$\lim_{\substack{\tau\rightarrow n\\ n\in\mathds Z}} I_K(\tau) = 1
\mbox{\hskip 0.2in and\hskip 0.2in}
I_K\left(\frac{m}{K}\right) = 0,\;\;\; 
m\in\mathds Z,\; m\mbox{ not a multiple of $K$.}$$
Now let $x(\tau) = x(\tau+T_0)$ be a periodic function of $\tau$
with period $T_0$ and set $x_j = x(jT_0/K)$ for $j = 0\ldots K-1$.
Then, defining
\begin{equation}
\wh{x}(\tau) := \sum_{j = 0}^{K-1} x_j I_K(\tau/T_0 - j/K),
\label{xi_def}
\end{equation}
we have, in the light of (i) and (ii) above, that $\wh{x}(kT_0/K)
= x_k = x(kT_0/K)$ for $k\in\mathds Z$. Hence, $\wh{x}(\tau)$
can be used to interpolate $x(\tau)$ given the values of $x(\tau)$
on a uniformly-spaced discrete set of values of $\tau$. 
In Appendix \ref{app:a} we show that the error in the
interpolation scheme described is $O\left({\rm e}^{-C_2 K}\right)$,
for some positive constant $C_2$.

In practice, for $\tau$ close to an integer, $I_K(\tau)$ is best
computed from a series expansion. Letting $\delta = \tau - [\tau]$,
with $[\tau]$ being the nearest integer to $\tau$, we then use
\begin{equation}
I_K(\delta) = 1 - \frac{1}{6}(K^2-1)\left[(\pi \delta)^2 - 
\frac{1}{60}(3K^2-7) (\pi \delta)^4 +
\frac{1}{2520}(3K^4 - 18K^2 + 31)(\pi \delta)^6\right] + O(\delta^8)
\label{I_series}
\end{equation}
whenever $\left|\delta\right| < \varepsilon_I$. Since the computations
are carried out to approximately 16 significant figures,
we allow a margin of safety by choosing $\varepsilon_I = 10^{-4}$.

The use of $I_K(\tau)$ for interpolation has other advantages,
amongst them that $\wh{x}(\tau)$ can be integrated in closed form,
and so, by implication, the
integral of $x(\tau)$ for all $\tau$ can be approximated. By defining
$$ J'_K(T) = \int_0^T{\mathrm d}\tau\,
\frac{\sin K\pi \tau}{\sin \pi \tau} $$
it is easy to show that
$$ J'_{K+2}(T) = J'_K(T) + 2\int_0^T {\mathrm d}\tau\,\cos (K+1)\pi
\tau = J'_K(T) + \frac{2}{(K+1)\pi}\sin (K+1)\pi T.$$
Now, since $K>0$ is odd and $J'_1(T) = T$, we have 
$$ J'_K(T) = T + \frac{1}{\pi}\sum_{i=1}^{(K-1)/2}
\frac{1}{i}\sin 2i\pi T.$$
Defining now $J_K(T) := \int_0^T {\mathrm d}\tau\, I_K(\tau) = J'_K(T)/K$,
we have
\begin{equation}
J_K(T) = \frac{T}{K} + \frac{1}{K\pi}\sum_{i=1}^{(K-1)/2}
\frac{1}{i}\sin 2i\pi T.
\label{J_K} \end{equation}
Next define $\wh{X}(T) = \int_0^T\mathrm{d}\tau\,\wh{x}(\tau)$. 
Integrating~(\ref{xi_def}) term-by-term, we obtain
\begin{equation}
\wh{X}(T) = T_0\sum_{i=0}^{K-1}x_i\left\{J_K\left(\frac{T}{T_0} -
\frac{i}{K}\right) + J_K\left(\frac{i}{K}\right)\right\} ,
\label{X_def} \end{equation}
where we have used the fact that $J_K(\tau)$ is an odd
function of $\tau$. In what follows, we therefore use $\wh{X}(T)$
to approximate $\int_0^T\mathrm{d}\tau\,x(\tau)$.

In a similar manner, defining
$E_K(\zeta, T) := \int_0^T\mathrm{d}\tau\,{\rm e}^{-\zeta\tau}I_K(\tau)$,
for constant $\zeta$, it can be shown that
\begin{equation}
E_K(\zeta, T) = \frac{1 - {\rm e}^{-\zeta T}}{\zeta K} +
\frac{2}{K}\sum_{i = 1}^{(K-1)/2} \frac{\zeta + {\rm e}^{-\zeta T}
\left( 2i\pi\sin 2i\pi T - \zeta\cos 2i\pi T\right)}
{\zeta^2 + 4\pi^2 i^2}.
\label{E_K}
\end{equation}
Hence, $\wh{X}_{\rm e}(\zeta, T) := \int_0^T\mathrm{d}\tau\,
{\rm e}^{-\zeta\tau}x(\tau)$ is given by
\begin{equation}
\wh{X}_{\rm e}(\zeta, T) = T_0\sum_{i=0}^{K-1} x_i\, {\rm e}^{-\zeta i/K}
\left\{E_K\left(\frac{T}{T_0} - \frac{i}{K}\right) -
E_K\left(-\frac{i}{K}\right)\right\}.
\label{IE_K}
\end{equation}

\subsection{Calculation of {\boldmath $a(\tau),\, b(\tau)$}}

Before we can compute $w_{11}(\tau)$, we need to find the
unperturbed limit cycle, its period, $T_0$, the periodic
function $\tilde{F}(\tau)$ and the mean of $F(\tau)$, $f_0$. 
These are all straightforward: we solve the unperturbed
equation~(\ref{eq:2.5}) numerically as described in
Section~\ref{nsol}, obtaining 
$U_0$, $T_0$ and the solution and its derivative over one period.
Since $u_0(\tau)$ is periodic, so is $h(u_0(\tau))$, 
and so we can use equation~(\ref{X_def}) to estimate $F(\tau)$
for any $\tau$. From $F(\tau)$ we can then obtain $f_0$,
and hence $\tilde{F}(\tau)$.

Computation of $w_{11}(\tau)$ can now be carried out from
equation~(\ref{eq:2.15}), but is complicated by the fact that,
for $\tau = i T_0/2$, $i\in\mathds Z$, the integrand is singular
and numerical integration techniques will break down.
Singularity in the integrand, which is periodic, also prevents 
us from using equation~(\ref{X_def}). To discuss this further,
let us define two subsets of $\mathds R$ as $S = \cup_{i\in\ZZZ} s_i$
with $s_i =  \left[i T_0/2 - r_c, i T_0/2 + r_c\right]$,
where $r_c\ll T_0/2$ is small and will be defined later; and $I =
\mathds R\setminus S$. We will then use a power series representation
for $w_{11}(\tau)$, $\tau\in S$, where the power series
converges `usefully' (the error term is less than the
maximum acceptable error) for $|\tau| \leq r_c$, with
Romberg integration~\cite{numrec}, a standard numerical
integration technique, being used for $\tau\in I$.

It should be pointed out here that we do not need to compute
$\bar{\tau}$  explicitly. Instead, we can set the lower limit
of the integral to any convenient value, $\tau_l$ say, provided we add
a suitable multiple of $\dot{u}_0(\tau)$; so our definition becomes
\begin{equation}
w_{11}(\tau) = \dot{u}_0(\tau) k_2 + 
\dot{u}_0(\tau)\int_{\tau_l}^\tau d\tau' c_{1}
\frac{{\rm e}^{-F(\tau')}}{\dot{u}_0^2(\tau')} ,
\label{w11def} \end{equation}
where the constant $k_2$, which depends on $\tau_l$,
will be chosen to ensure continuity.

In practice, we only need to know $w_{11}(\tau)$ over a length of
time consisting of two periods of $u_0(\tau)$,
and so we calculate it for $\tau\in [0, 2 T_0]$: from
Appendix A in~\cite{BDG3}, we know that $w_{11}(\tau)$ is well-defined
even at $\tau = 0$, which we take to be our value of $\tau_l$.

We derive the formal power series for $w_{11}(\tau)$ by using the
method of Frobenius to solve the differential equation~(\ref{eq:2.5}), 
with initial conditions chosen so as to ensure that the solution is
on the limit cycle. Thus, $u_0(0) = U_0$, $\dot{u}_0(0) = 0$, from which we
obtain the power series in $\tau$ for $u_0(\tau)$ and hence, using
term-by-term differentiation, for $\dot{u}_0(\tau)$. The latter is
$$ \dot{u}_0(\tau) = U_0\left(\alpha -\beta + \beta U_0^2\right)
\left[-\tau + \left(1- \beta + 3\beta U_0^2\right) 
\frac{\tau^2}{2}\right] + O\left(\tau^3\right). $$
By Remark~1 in~\cite{BDG3}, $c_1 = -r_1$, which is the coefficient
of $\tau$ in the above series for $\dot{u}_0(\tau)$, and so $c_1 =
U_0(\alpha - \beta + \beta U_0^2)$. Using the series for $u_0(\tau)
= U_0 + \int_0^\tau d\tau' \dot{u}_0(\tau')$, and term-by-term
integration, we can also find power series for $F(\tau)$,
${\rm e}^{-F(\tau)}$ and hence for the integrand
${\rm e}^{-F(\tau)}/\dot{u}_0^2(\tau)$. Integrating
this series from 0 to $\tau$ term-by-term, and multiplying by
the series for $c_1 \dot{u}_0(\tau)$, we obtain $w_{11}(\tau) =
1 - \left(1- \beta + 3\beta U_0^2\right) \tau/2 +
\left(1-2\alpha + \beta^2(1-3U_0^2)^2\right) \tau^2/4 + O(\tau^3)$.
Finally, we apply the remaining condition on $\dot{w}_{11}$,
that is, $\dot{w}_{11}(0) = 0$,  which forces the choice of $k_2$ 
in equation~(\ref{w11def}) to be such that $-k_2 U_0
\left(\alpha -\beta + \beta U_0^2\right) - \left(1- \beta + 3\beta
U_0^2\right)/2 = 0$. This gives
\begin{equation}
w_{11}^{\mbox{\scriptsize{ser}}}(\tau) 
\approx 1 + \sum_{j = 2}^M R_j \tau^j + O\left(\tau^{M+1}\right) ,
\label{w11ser} \end{equation}
where $2R_2 = \alpha - \beta + 3\beta u_0^2$, $6 R_3 = \alpha -
\beta(1 + \alpha - \beta) + 3\beta(1 + 3\alpha - 4\beta)u_0^2 +
15\beta^2u_0^4$ and so on. This is a truncation of the series
actually used for $\tau\in s_0$. Using computer
algebra, it is possible to extend this series to at least the
term of order $\tau^{10}$, expressing each coefficient of $\tau$
as a polynomial in $\alpha, \beta$ and $U_0$ --- that is,
without assuming numerical values for these parameters --- although
the higher order coefficients become quite complicated.

The series~(\ref{w11ser}) can be used to estimate
$w_{11}(\tau)$ for $\tau\in s_j$, $j > 0$, provided (i) the
value given by the series is multiplied by $(-1)^j
{\rm e}^{-j f_0 T_0/2}$ and (ii) $k_2$ is chosen so as to ensure
continuity across the boundary of $s_j$.
The term $(-1)^j$ in the correction factor arises as a result 
of the property of $u_0(\tau)$ in equation~(\ref{eq:2.6}),
and the exponential factor comes
from the definition of $w_{11}(\tau)$, equation~(\ref{eq:2.15}). 
Hence, $w_{11}(\tau)$ is estimated as 
\begin{equation}
w_{11}(\tau) = k_2 \dot{u}_0(\tau) + (-1)^j {\rm e}^{-j f_0 T_0/2} 
w_{11}^{\mbox{\scriptsize{ser}}}(\tau - j T_0/2) + O\left(\tau^{M+1}\right)
\label{w11est} \end{equation}
for $\tau\in s_j$.

\begin{figure}[htbp]
\centering
\includegraphics*[angle=0,width=5.9in]{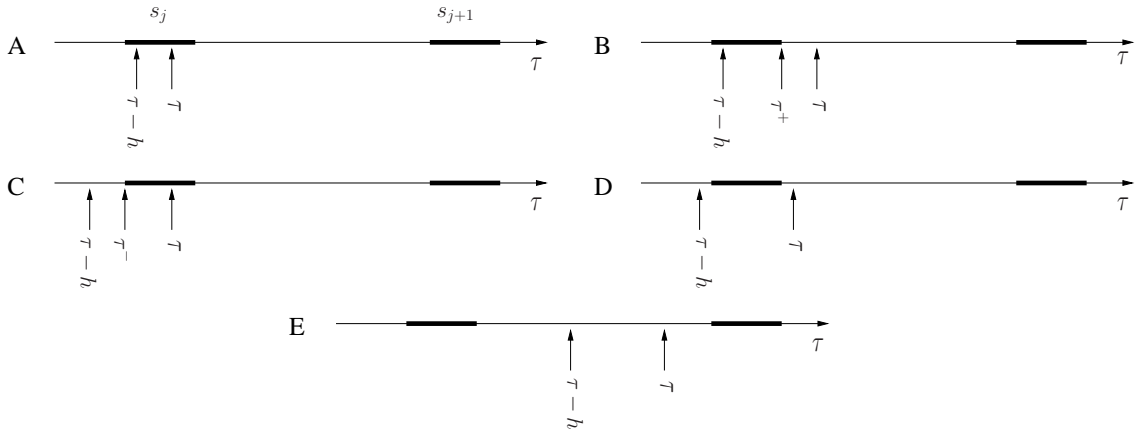}
\caption{The cases A--E to be considered when calculating $w_{11}(\tau)$.
The subsets $s_j = [j T_0/2 - r_c, j T_0/2 + r_c]$ and $s_{j+1}$
are shown as thick line segments. It is assumed that $w_{11}$
is being computed at equally spaced time steps of width $h$,
that $w_{11}(\tau-h)$ has just been calculated, and that $w_{11}(\tau)$ is
now to be found. How this is done depends on the
relationship of $\tau$ to $\tau-h$ --- see text.}
\label{fig:cases}
\end{figure}

In order to compute $b(\tau)$, we need to know $w_{11}(\tau)$
for $\tau\in [0, 2T_0]$ --- see equation~(\ref{bcalc}) --- and hence
we calculate $w_{11}(\tau)$ at equally spaced points $0, h, 2h,
\ldots, 2 T_0$, in that order, where $2T_0/h$ is an integer. The point 
$\tau=0$ is in $s_0$ and so the truncated series is used here
(with $k_2 = 0$). For $\tau>0$, various different cases exist,
and these are illustrated in figure~\ref{fig:cases},
in which $\tau$ is the time at which $w_{11}$ is to be calculated,
and it is assumed that it has already been calculated at $\tau-h$.

\begin{itemize}
\item In case A, $\tau\in s_j$, so the series is used, 
with the current value of $k_2$.
\item In case B, first $w_{11}(\tau^+)/\dot{u}_0(\tau^+)$, 
where $\tau^+ = jT_0/2 + r_c$,
is calculated from the truncated series. To this is added a 
numerical estimate of $\int_{\tau^+}^\tau {\rm d}\tau'
c_1 {\rm e}^{-F(\tau')}/\dot{u}_0^2(\tau')$, and the result
multiplied by $\dot{u}_0(\tau)$ to obtain an estimate of $w_{11}(\tau)$.
\item In case C, roughly the reverse happens. Define $\tau^- =
jT_0/2 - r_c$. Then numerical integration is used to
estimate $w_{11}(\tau^-)$, from which $k_2$ can be found, by assuming
continuity across the boundary $\tau = \tau^-$. Since the
appropriate value of $k_2$ is now known, the truncated series
can be used to estimate $w_{11}(\tau)$.
\item In case D, compute as in C, followed by B.
\item In case E, straightforward numerical integration alone is used.
\end{itemize}

In this way, $w_{11}(ih)$ is computed for $i = 0, 1, 2, \ldots, 2T_0/h$,
and it is now a simple matter to extract $a(\tau)$
and $b(\tau)$ at the points $\tau = ih$, $i = 0, 1,\ldots, T_0/h$,
so that their values at any $\tau$ can be found
by interpolation. From equation~(\ref{eq:2.18}), we have that
$w_{11}(\tau) = a(\tau) + {\rm e}^{-f_0 \tau} b(\tau)$ 
and $a(\tau) = \gamma \dot{u}_0(\tau)$. Since
$a(\tau)$ and $b(\tau)$ both have period $T_0$, we have
\begin{equation}
b(\tau) = {\rm e}^{f_0 \tau} \frac{w_{11}(\tau) -
w_{11}(\tau + T_0)}{1 - {\rm e}^{-f_0 T_0}}
\label{bcalc} \end{equation}
and, knowing $w_{11}(\tau)$ for $\tau\in[0, 2T_0]$,
we can now compute $b(\tau)$ for $\tau\in[0, T_0]$.
Having computed $b(\tau)$, we can use for instance
the method of least squares to estimate the value of $\gamma$:
that is, we find the value of $\gamma,
\hat{\gamma}$, that minimises
$$ \sum_{i=0}^{T_0/h}\left[w_{11}(ih) - {\rm e}^{-f_0 ih} b(ih)-
\hat\gamma\dot{u}_0(ih)\right]^2,$$
which is
\begin{equation}
\hat\gamma = \frac{\sum\dot{u}_0(ih)
\left[w_{11}(ih) - {\rm e}^{-f_0 ih} b(ih)\right]}
{\sum\dot{u}_0^2(ih)},
\label{eq:gamma} \end{equation}
where the sums go from $i = 0$ to $T_0/h$.
This completes the calculation of $a(\tau)$ and $b(\tau)$.

\subsection{Illustrative results}
\label{sec:3.3}

Illustrative results are now given for the case $\alpha = 5$,
$\beta = 4$ and $f(\tau) = \sin\tau$. All the computations were carried out 
using double precision arithmetic. For interpolation, $K = 151$ equally
spaced points were used; the series for $I_K(\delta)$
was used if $|\delta| < \varepsilon_I = 10^{-4}$.
In series~(\ref{w11ser}), $M = 10$. In the definition of $s_j$, $r_c =
10^{-2}$, and the fractional accuracy chosen for Romberg integration
was $10^{-12}$.  With these parameters, and $\rho = 2$, we find $T_0 \approx
3.698939867513906$, $U_0 \approx 0.979106186033891$,
$f_0 \approx 0.757499334158$ and $\hat\gamma = -54.855909271256$.
Having computed $a(\tau)$ and $b(\tau)$, we can then estimate $A$, 
first of all from equation~(\ref{eq:2.24}), 
using Romberg integration: this gives $A = 16.0813516305191$.
Using equation~(\ref{X_def}) to carry out the integration, we obtain
$A = 16.0813516305189$.  Furthermore, we have from 
equation~(\ref{eq:2.28}) that $A = -r_1\rho\Omega_0 = U_0\left(\alpha -
\beta + \beta U_0^2\right)/(\rho\Omega_0) = 16.0813516307791$.
These estimates agree with each other to 11 significant figures,
thereby verifying the numerical techniques used to obtain them.

The calculation of $B_{11}\ldots B_{32}$ and $D_{1}$, $D_{2}$, defined after
equation~(\ref{eq:2.27}), now follows
straightforwardly from equations~(\ref{eq:2.25}). The only
point to note is that these integrals can be zero, which gives
problems in the error control scheme used for numerical integration.
To overcome this, the integration is done in two parts,
from 0 to $2\pi p z$ and from $2\pi p z$ to $2\pi p$,
with $z$ approximately, but not exactly, one half.
We then find, for the above parameters, that,
when $p = 1$, $D_1, D_2 \approx 10^{-12}$. On the other hand,
with $p = 2$, we find $D_1 \approx 8.11989\times 10^{-2}$ and 
$D_2 \approx -5.20174\times 10^{-1}$; for $p = 4$,
$D_1 = -3.79022\times 10^{-2}$, $D_2 = 2.74434\times 10 ^{-1}$.

\zerarcounters
\section{Numerical results}
\label{sec:5}

To extend our analytical results to large values of $\mu$ we compute
a set of \Arnold\ tongues numerically. We make several comparisons
between the theoretical predictions and the computational results,
which provide information about the computational accuracy and
the range of validity of some theoretical estimates.

\subsection{\Arnold\ tongues}
\label{sec:5.1}

We computed the \Arnold\ tongues of system~(\ref{eq:2.1}) for
the 15 strongest resonances using the algorithms from~\cite{FS}
with two types of forcing: Harmonic forcing with
\begin{equation}
f(t) = \sin t
\label{eq:num:harmf} \end{equation}
as considered in~\cite{BDG3}, and the forcing function
\begin{equation}
f(\tau) = \frac{\lambda^2-1}{\lambda}
\sum_{k = 1}^{\infty} \frac{\sin k\tau}{\lambda^{k}} 
= \frac{ \left( \lambda^{2} - 1 \right) \sin \tau}{\lambda^{2} +
1 - 2 \lambda \cos \tau} , \quad \lambda > 1,
\label{eq:num:genf} \end{equation}
containing all harmonics. Note that $f$ in~(\ref{eq:num:genf}) is smooth
and $f(\tau)\in[-1, 1]$; hence, it can be used as a direct replacement
for the sine function in~(\ref{eq:num:harmf}). The relative
strength of the harmonics can be adjusted by varying $\lambda$,
since one has $\hat f_{\nu}=\Phi(\lambda)\,\lambda^{-|\nu|}$,
with $\Phi(\lambda)=(\lambda^{2}-1)/2{\rm i}\lambda$.

\begin{figure}
\centering
\includegraphics[width=0.85\linewidth]{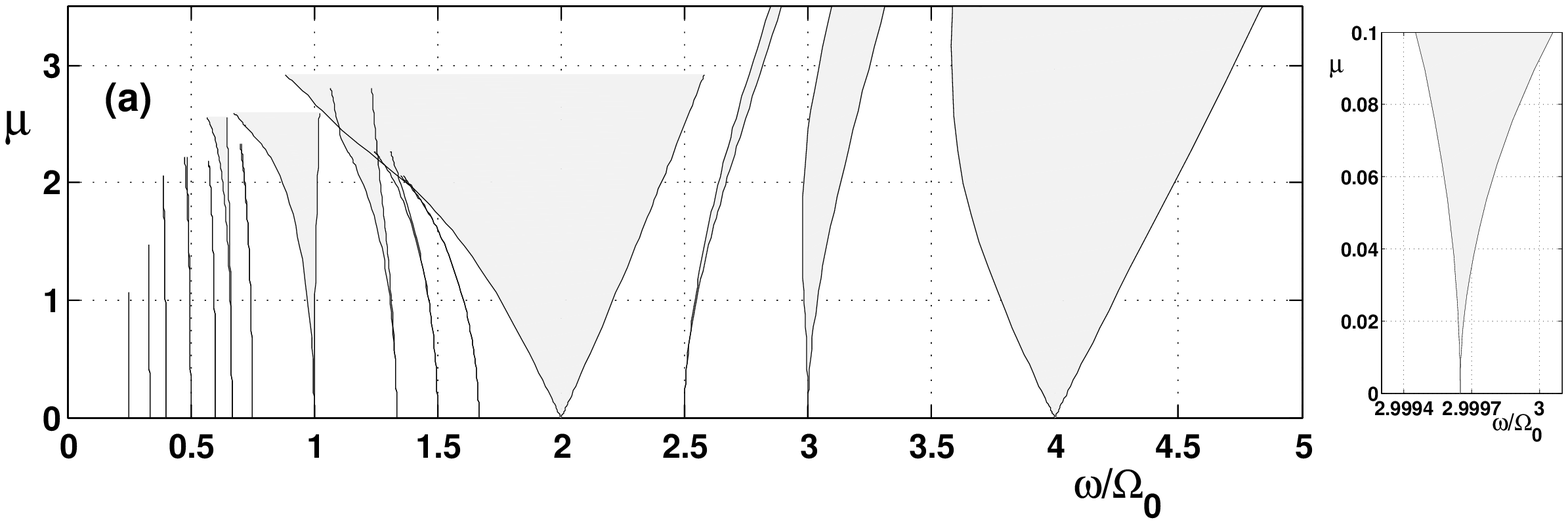} \\[1em]
\includegraphics[width=0.85\linewidth]{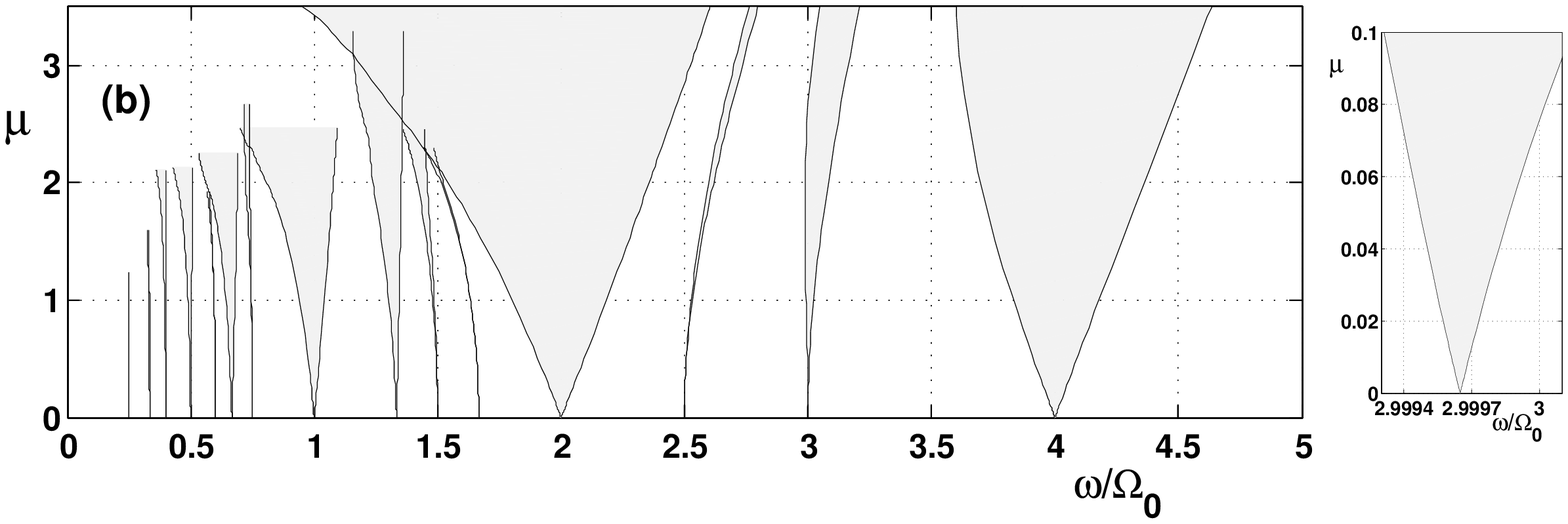}
\caption{Some \Arnold\ tongues in the $(\omega/\Omega_{0},\mu)$-plane for
the ILFD for $\alpha=5$ and $\beta=4$ with (a) $f(\tau)=\sin\tau$
and (b) $f(\tau)$ given by~(\ref{eq:num:genf}) with $\lambda=2$.
For these parameter values we have $\Omega_{0}=1.69864489$,
$T_{0}=2\pi/\Omega_{0}=3.69893987$, $A\rho=2.78668166$,
$D_{1}=0.00703534$ and $D_{2}=-0.0450695$.}
\label{fig:tongues}
\end{figure}

The results of both computations for the parameter values $\alpha=5$
and $\beta=4$ are illustrated and compared in figure~\ref{fig:tongues}.
Note that in this case we have $\Omega_{0}=1.698645\dots$ and, hence,
$\om_{\rm c}(2)=3.397290\dots$ and $\om_{\rm c}(4)=6.794580\dots$
As explained in detail in~\cite{FS} these tongues are computed
by continuation of so-called constant-$\mu$ cross sections, starting
at the tips. To facilitate our subsequent computations of the order
of contact and opening angles we started with an extremely small
continuation step size to obtain a large number of points very close
to the tips for later fitting. For each tongue we performed 150
continuation steps. The computation of most tongues terminated
by either reaching the computational boundary $\mu=3.5$ or by
exceeding the maximal number of 150 continuation steps. However,
the computation of some tongues, most notably of the 2:1 tongue,
seems to end due to limitations of the algorithm we use;
see~\cite{FS} for more details. We did not pursue a further
investigation, because we are mainly interested in the size
and location of the 2:1 and 4:1 tongues for moderate $\mu$
and in investigating the behaviour at the tips of all tongues,
for which we obtained sufficient data.

The bifurcation diagrams shown in figure~\ref{fig:tongues} are clearly
dominated by the strongest resonances occurring for $\rho=2$ and $\rho=4$.
A continuation of the frequency-locked sub-harmonic solutions along the
centre-lines $\omega/\Omega_0=2$ and $\omega/\Omega_0=4$ inside
these tongues revealed that these solutions remain attracting
for $\mu\leq 3.5$ at least along these centre-lines. On the other hand,
we observe at the left-hand boundary of the 2:1 tongue that this tongue
overlaps with other tongues. Hence, in these overlapping regions we
might find multi-stability. To the right-hand side of the 2:1 tongue
no such phenomenon is apparent in these figures.

The small plots to the right of figures~\ref{fig:tongues}~(a) and~(b)
show enlargements of the tip of the 3:1 tongue illustrating the effect
of forcing with and without all harmonics present as predicted in 
Section~\ref{sec:4.1}. In figure~\ref{fig:tongues}~(a) we observe a 
high-order (quadratic) contact of the boundaries of the tongue, 
while in figure~\ref{fig:tongues}~(b) the two boundaries intersect 
transversally. Note that the slight shift to the left of 
$\omega/\Omega_0=3$ is due to the discretisation error of the 
periodic solutions; see~\cite{FS} for technical details. It is 
remarkable, however, that our computations accurately capture the 
predicted high-order behaviour despite this approximation error, 
which is orders of magnitudes larger than the width of most tongues 
at their tips.

\begin{figure}
\centering\small
\begin{tabular}{cc}
\includegraphics[width=0.45\linewidth]{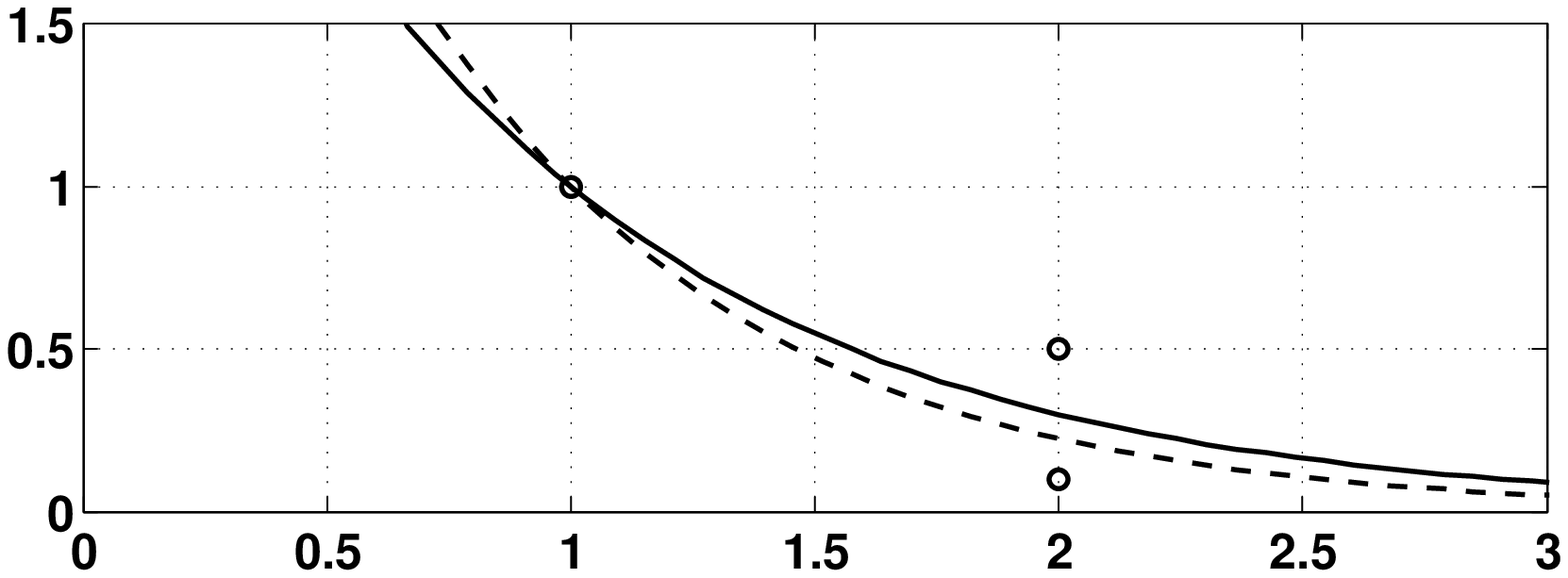} &
\includegraphics[width=0.45\linewidth]{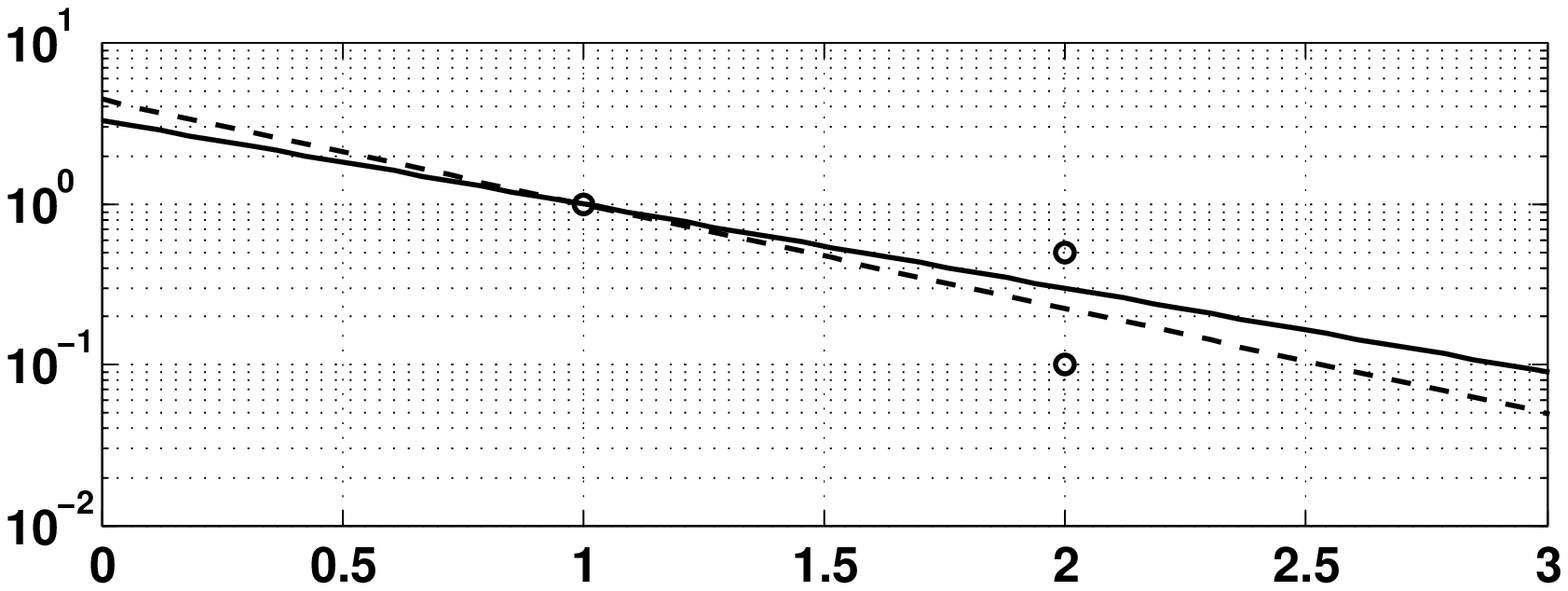} \\[1ex]
\bf (a) & \bf (b)
\end{tabular}
\caption{Comparison of linear (dashed) and non-linear (solid)
fits for a simple test example.}
\label{fig:comp_fit}
\end{figure}

To quantify our findings and for comparison with the analytical 
predictions we developed a simple adaptive non-linear fitting 
algorithm for the width-function $\Delta\omega(\rho)$ to a 
monomial $\Delta\omega(\rho)=a\mu^b$ with $a$ and $b$ unknown. 
One might argue that one could use a linear fit to logarithmic 
data of the form $\ln(\Delta\omega(\rho))=\ln a + b\ln\mu$ 
to compute estimates for $a$ and $b$. However, this leads to
biased estimates as we illustrate in figure~\ref{fig:comp_fit}, 
where we compare the results of a non-linear fit (solid) with a 
linear fit (dashed) to the function $y=a\exp(bx)$. Only the 
non-linear fit is a useful fit to the data 
as figure~\ref{fig:comp_fit}~(a) clearly illustrates, the linear fit 
here being biased towards lower function values. There are 
several reasons why a linear fit to logarithmic data is 
inappropriate, the most important ones being that the two 
least-squares residual functions $\|Y-a\exp(bX)\|_2^2$ 
and $\|\ln Y-(\ln a + bX)\|_2^2$ have different minimisers, 
and that $Y$ and $\ln Y$ do not have the same 
error distribution. Another suggestion could be to compute 
a linear fit to a polynomial $p_n(\mu) = a_0+a_1\mu+\dots +a_n\mu^n$ 
of sufficiently high order~$n$. However, we found that this leads 
to a least-squares problem that is so ill-conditioned that 
round-off errors become amplified to order one, that is, 
the fitted coefficients are essentially meaningless. 
A way out is to use orthonormal polynomials
as base functions instead of monomials of the form $a_k\mu^k$.
However, since our non-linear fit worked sufficiently well we did 
not pursue this further.

For our computations we used the weighted least-squares residual function
\begin{equation*}
F(a,b) := \left\|W(\mu)\left(\Delta\omega(\rho)-
a\mu^{b}\right)\right\|_{2}^{2},
\end{equation*}
where we used the weight-function $W(\mu)=1/\mu$ to penalise errors 
closer to the tip $\mu=0$. The numerical data is normalised 
to the unit square, that is, we compute a fit to 
$\Delta\omega(\rho)/\max\{\Delta\omega(\rho)\}$ and $\mu/\max\{\mu\}$. 
We computed estimates for $a$ and $b$ by applying Newton's method 
to the equations $\partial F(a,b)/\partial(a,b)=0$ and rescaled 
the computed coefficients to fit the original data. The size of the 
fitting interval $\mu\in[0,\mu_{\mathrm{fit}}]$ was computed adaptively. 
We started with an initial fitting interval $\mu_{\mathrm{fit}} = 
\max\{0.1, \mathrm{argmax} \{\Delta\omega(\rho) < 10^{-4} \} \}$, 
that is, we used either $\mu_{\mathrm{fit}}=0.1$ or the largest value 
of $\mu$ such that the width of the tongue was less than $10^{-4}$.
 The fitting interval was accepted if the least squares error $F(a,b)$ 
for the normalised data was less than $10^{-3}$ and reduced successively 
if the error was larger. Within the fitting interval we excluded points
for which $\Delta\omega(\rho)$ was zero within numerical accuracy.
Table~\ref{tab:1} summarises our results 
for both forms of forcing. Each row states the fitted monomial 
representing the leading order term together with the fitting interval 
and the number $N_{\mathrm{fit}}$ of data points this monomial 
was fitted to. These computations agree extremely well with the 
theoretical predictions and also verify that our fitting algorithm 
is suitable to capture the leading-order behaviour accurately.

\begin{table}
\centering\small
\setlength\tabcolsep{5pt}
\begin{tabular}{|l|llr|llr|}
\hline &&& &&& \\[-1em]
& \multicolumn{3}{c|}{$f(\tau)=\sin(\tau)$} &
	\multicolumn{3}{c|}{$f(\tau)=
	\frac{(\lambda^{2} - 1 ) \sin \tau}
	{\lambda^{2} + 1 - 2 \lambda \cos \tau}$, $\lambda=2$} \\[0.4ex]
\cline{2-7} &&& &&& \\[-1em]
$p\!:\!q$ &
	\multicolumn{1}{|c}{$\Delta\omega(p/q)$} &
		\multicolumn{1}{c}{$\mu_{\mathrm{fit}}$} &
		$N_{\mathrm{fit}}$ &
	\multicolumn{1}{|c}{$\Delta\omega(p/q)$} &
		\multicolumn{1}{c}{$\mu_{\mathrm{fit}}$} &
		$N_{\mathrm{fit}}$ \\[0.3ex]
\hline 
\hline &&& &&& \\[-1em]
1:4 &
$2.770 \times {10}^{-7}\: \mu^{8.076}$ & $8.763 \times {10}^{-1}$ & 48 &
$5.571 \times {10}^{-4}\: \mu^{1.001}$ & $4.481 \times {10}^{-3}$ & 72 \\
1:3 &
$1.055 \times {10}^{-5}\: \mu^{6.052}$ & $8.971 \times {10}^{-1}$ & 40 &
$2.971 \times {10}^{-3}\: \mu^{1.001}$ & $5.932 \times {10}^{-3}$ & 73 \\
2:5 &
$7.554 \times {10}^{-5}\: \mu^{5.038}$ & $8.294 \times {10}^{-1}$ & 36 &
$7.128 \times {10}^{-3}\: \mu^{1.001}$ & $7.913 \times {10}^{-3}$ & 73 \\
1:2 &
$5.188 \times {10}^{-4}\: \mu^{4.014}$ & $6.529 \times {10}^{-1}$ & 40 &
$1.786 \times {10}^{-2}\: \mu^{1.001}$ & $1.068 \times {10}^{-2}$ & 74 \\
3:5 &
$3.093 \times {10}^{-8}\: \mu^{10.19}$ & $8.839 \times {10}^{-1}$ & 14 &
$1.057 \times {10}^{-5}\: \mu^{1.001}$ & $7.612 \times {10}^{-4}$ & 58 \\
\hline &&& &&& \\[-1em]
2:3 &
$4.096 \times {10}^{-3}\: \mu^{3.003}$ & $2.662 \times {10}^{-1}$ & 49 &
$4.780 \times {10}^{-2}\: \mu^{1.001}$ & $1.612 \times {10}^{-2}$ & 74 \\
3:4 &
$2.409 \times {10}^{-6}\: \mu^{8.016}$ & $7.462 \times {10}^{-1}$ & 16 &
$5.101 \times {10}^{-5}\: \mu^{0.9990}$ & $7.016 \times {10}^{-4}$ & 55 \\
1:1 &
$4.828 \times {10}^{-2}\: \mu^{2.000}$ & $9.099 \times {10}^{-2}$ & 76 &
$1.467 \times {10}^{-1}\: \mu^{1.001}$ & $3.179 \times {10}^{-2}$ & 76 \\
4:3 &
$9.979 \times {10}^{-3}\: \mu^{3.001}$ & $1.540 \times {10}^{-1}$ & 47 &
$5.016 \times {10}^{-2}\: \mu^{1.001}$ & $1.097 \times {10}^{-2}$ & 70 \\
3:2 &
$1.112 \times {10}^{-3}\: \mu^{4.014}$ & $5.338 \times {10}^{-1}$ & 38 &
$1.688 \times {10}^{-3}\: \mu^{1.001}$ & $1.973 \times {10}^{-3}$ & 61 \\
\hline &&& &&& \\[-1em]
5:3 &
$5.049 \times {10}^{-5}\: \mu^{5.971}$ & $2.588 \times {10}^{-1}$ & 16 &
$2.854 \times {10}^{-5}\: \mu^{0.9990}$ & $5.950 \times {10}^{-4}$ & 55 \\
2:1 &
$7.556 \times {10}^{-1}\: \mu^{1.000}$ & $9.684 \times {10}^{-2}$ & 80 &
$6.280 \times {10}^{-1}\: \mu^{1.000}$ & $9.024 \times {10}^{-2}$ & 79 \\
5:2 &
$1.595 \times {10}^{-2}\: \mu^{3.971}$ & $2.039 \times {10}^{-1}$ & 35 &
$1.832 \times {10}^{-4}\: \mu^{0.9990}$ & $4.539 \times {10}^{-4}$ & 50 \\
3:1 &
$1.024 \times {10}^{-1}\: \mu^{1.999}$ & $8.960 \times {10}^{-2}$ & 77 &
$1.331 \times {10}^{-2}\: \mu^{0.9994}$ & $2.410 \times {10}^{-2}$ & 73 \\
4:1 &
$7.957 \times {10}^{-1}\: \mu^{0.9999}$ & $9.719 \times {10}^{-2}$ & 80 &
$5.968 \times {10}^{-1}\: \mu^{0.9999}$ & $9.757 \times {10}^{-2}$ & 79 \\
\hline
\end{tabular}
\caption{Leading contributions to the plateau widths corresponding
to the main resonances as they appear in figure~\ref{fig:tongues} 
from left to right. These coefficients were obtained by fitting 
the monomial $a\mu^b$ to the numerically computed values 
for~$\Delta\omega(p/q)$ over the interval $\mu\in[0,
\mu_{\mathrm{fit}}]$ on~$N_{\mathrm{fit}}$ data points.}
\label{tab:1}
\end{table}

First of all, our computations of the width of the locking intervals 
are in alignment with our theoretical results
\begin{equation*}
\Delta\omega(p/q) = \left\{\begin{array}{ll}
O(\mu^{k}) & \mathrm{for~} p \mathrm{~even}, \\
O(\mu^{2k}) & \mathrm{for~} p \mathrm{~odd},
\end{array}\right.
\end{equation*}
for harmonic forcing~(\ref{eq:num:harmf}), and with
\begin{equation*}
\Delta\omega(p/q) = O(\mu)
\end{equation*}
for general forcing~(\ref{eq:num:genf}) containing all harmonics. 
Furthermore, for the main \Arnold\ tongues corresponding to 
the resonances 2:1 and 4:1, equation~(\ref{eq:2.27}) valid for
harmonic forcing reduces to
\begin{equation}
\gotD_{1}(\tau_{0}) = \frac{1}{A} \left( D_{1} \cos \tau_{0} +
D_{2} \sin \tau_{0} \right) ,
\label{eq:5.3} \end{equation}
where, for the 2:1 resonance, $A=16.0814$, $D_{1}=8.11989 \times 10^{-2}$,
and $D_{2}= -5.20174 \times 10^{-1}$; see also Section~\ref{sec:3.3}.
An easy computation gives
\begin{equation}
\max_{0 \le \tau_{0} \le 2\pi} \gotD_{1}(\tau_{0}) =
- \min_{0 \le \tau_{0} \le 2\pi} \gotD_{1}(\tau_{0}) = M ,
\qquad M = 0.0327381 .
\label{eq:5.4} \end{equation}
For example, for the boundaries of the tongue corresponding
to the 2:1 resonance we find from table~\ref{tab:1} 
that $(\tan\theta_{1}(2) + \tan\theta_{2}(2))/2 \approx 
0.7556/2 = 0.3778$ in agreement with~(\ref{eq:2.32}),
which gives $\rho^{2}\Omega_{0}^{2} M=0.37785$. Similarly, 
for the 4:1 resonance we have $A = 32.1627$, $D_{1} =
-3.79022\times 10^{-2}$,
$D_{2} = 2.74434\times 10^{-1}$, giving $M = 8.6137\times 10^{-3}$. 
We compute that $\rho^{2}\Omega_{0}^{2} M= 0.39766$, again in 
agreement with the result
in table~\ref{tab:1} that $(\tan\theta_{1}(4) + \tan\theta_{2}(4))/2 
= 0.7959/2 = 0.39795$.

Also for the secondary resonances $p\!:\!1$, with $p$ odd,
the agreement between the numerical results and the
analytical predictions is satisfactory. The second order
computation, performed according to the analysis in Section \ref{sec:4.2},
gives, for the 1:1 and 3:1 resonances, the values
$\Delta(1)\approx 4.8246 \times 10^{-2}$ and
$\Delta(3)\approx 1.0269 \times 10^{-1}$, to be compared
with the values $4.828 \times 10^{-2}$ and $1.024 \times 10^{-1}$
in table~\ref{tab:1}.

\begin{figure}[htbp]
\centering
\subfigure[]{
    	\includegraphics*[angle=0,width=2.5in]{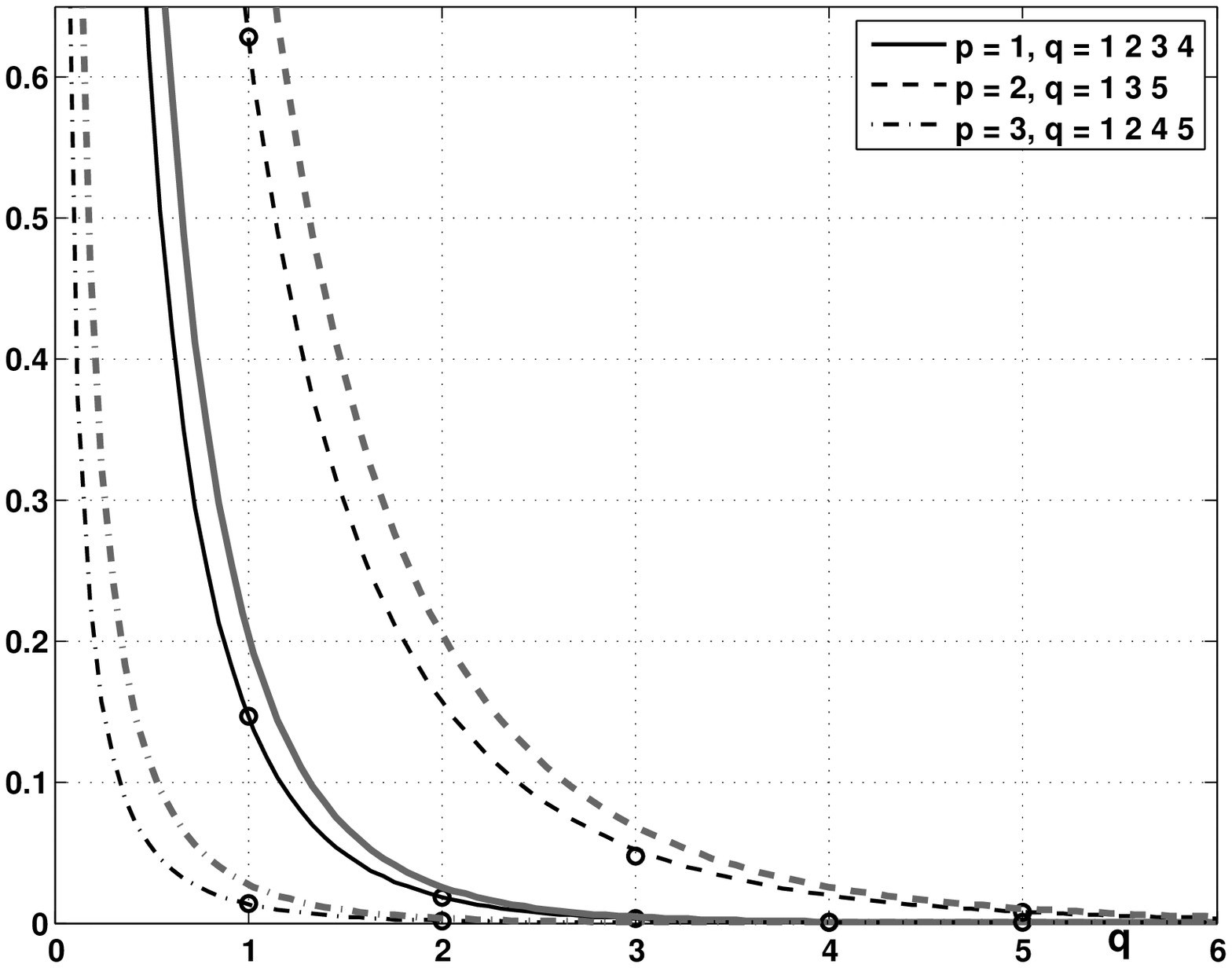}}
\hspace{0.5in}
\subfigure[]{
    	\includegraphics*[angle=0,width=2.5in]{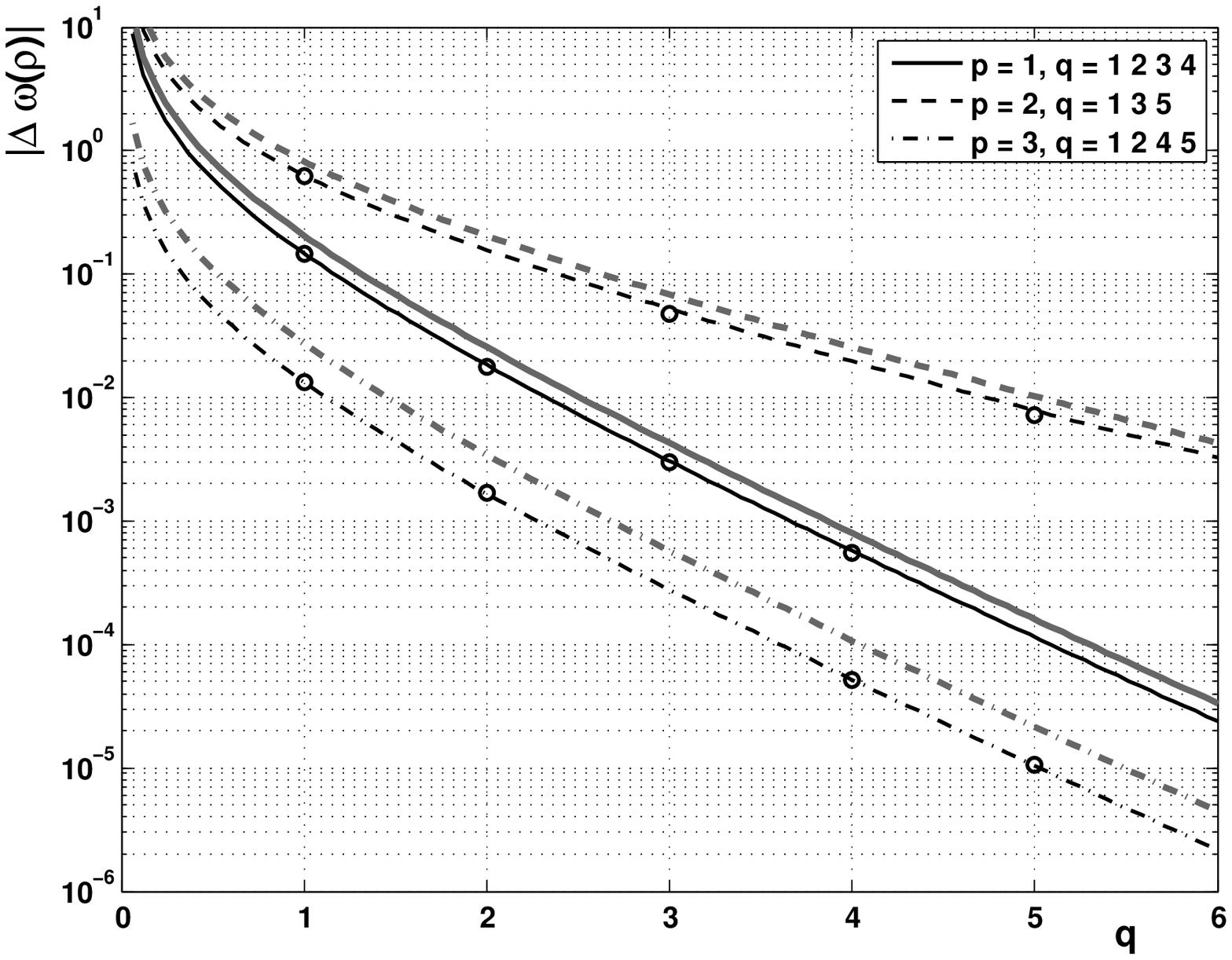}}
\vspace{-0.3cm}
\caption{Tongue widths $|\Delta\omega(\rho)|$ for fixed $p$ and
varying $q$. The black curves are plots of~(\ref{eq:5.5}) using the
constants $c$ from the numerical data and the grey curves are plots
of~(\ref{eq:5.5}) using the constants $c$ predicted by the theory.
The discrepancies are due to the small values of $q$ are are
expected to be asymptotically zero for large $q$.}
\label{fig:fit_exp}
\end{figure}

For the all-harmonics forcing~(\ref{eq:num:genf}), 
the plateau widths as given in table~\ref{tab:1} are consistent 
with the scaling law~(\ref{eq:4.10}). If we fix $p$ to be an odd integer,
then $\nu_{0}\rho=\nu_{0}p/q=2p$, hence $\nu_{0}=2q$,
and we find $|\hat f_{\nu_{0}}|\,|\overline K_{\nu_{0}\rho}(\rho)|
= \Phi(\lambda) \,\lambda^{-2q} |\overline K_{2p}(p/q)|$. If, on the 
contrary, we fix $p$ to be an even integer then $\nu_{0}\rho=\nu_{0}p/q=p$, 
hence $\nu_{0}=q$, so that we obtain $|\hat f_{\nu_{0}}|\,|\overline
K_{\nu_{0}\rho}(\rho)| = \Phi(\lambda)\,\lambda^{-q} |\overline K_{p}(p/q)|$.
When inserted into~(\ref{eq:4.10}), this leads to
\begin{equation}
|\Delta\omega(\rho)| \approx \left\{\begin{array}{ll}
c \mu / (q2^{2q}) & \mbox{for~} p \mbox{~odd},\\[1ex]
c \mu / (q2^q) & \mbox{for~} p \mbox{~even},
\end{array}\right.
\label{eq:5.5}
\end{equation}
with the constant $c=c(p)$ independent of $q$. The constant $c$ can be
computed  using the theory. A comparison between~(\ref{eq:5.5})
and~(\ref{eq:4.10}) gives $c(p)=\ln 2\Omega_{0}|\bar r_{1}|^{-1}
\Phi(\lambda)\,p\,|\overline K_{\nu_{0}\rho}(p/q)|$,
with $\nu_{0}\rho=2q$ for odd $p$ and $\nu_{0}\rho=q$ for even $p$.
For $p=1, 2, 3$ we compute $c(p)=0.82$, $1.64$, $0.11$, respectively.
These estimates are consistent with our numerical data
in table~\ref{tab:1}. Fitting our data to function~(\ref{eq:5.5})
we obtain the numerical estimates $c(p)=0.5867$, $1.255$, $0.05326$,
respectively, which is in good agreement considering 
the limited numerical accuracy and that~(\ref{eq:5.5}) is valid for 
$q\to\infty$. Figure~\ref{fig:fit_exp} shows a comparison between the 
theoretically and numerically obtained width functions.

\subsection{Width of plateaux as a function of $\boldsymbol\alpha$
and $\boldsymbol\beta$}
\label{sec:5.2}

\begin{figure}[htbp]
\centering
\includegraphics*[angle=0,width=2.7in]{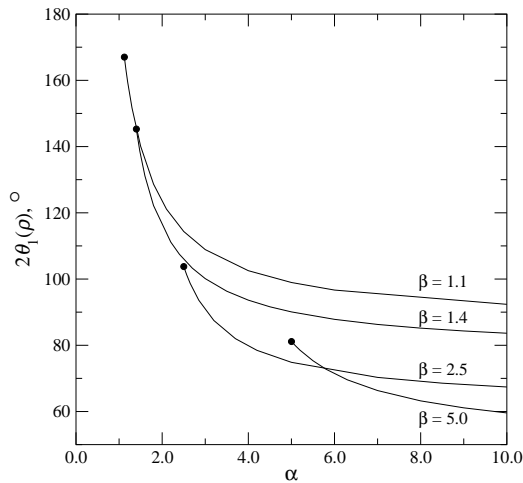}
\caption{A plot of $2\theta_1(2)$, the opening angle (degrees) of the
\Arnold\ tongue for $\rho = 2$, against $\alpha$, for various values
of $\beta$, with $\alpha > \beta > 1$. The angles have been computed
numerically from equation~(\ref{eq:2.36}).}
\label{fig:tanth}
\end{figure}

Of practical importance is the rate at which the width of a
given locking intereval increases with increasing $\mu$.
The locking regions, i.e. the \Arnold\ tongues,
are cone-shaped and the vertical angle of the cone,
$2\theta_1(\rho)$, which is a measure of width growth rate,
depends on the parameters $\alpha$ and $\beta$. This angle can
be computed from equation~(\ref{eq:2.36}) for a given $\rho$, and
typical results are given in figure~\ref{fig:tanth} for $\rho = 2$ and
a variety of values of $\beta$, with $\alpha\in(\beta, 10]$.

\zerarcounters
\section{Conclusions}
\label{sec:6}

In this paper we have investigated both analytically and numerically
the  structure of the \Arnold\ tongues for a resonant injection-locked
frequency divider (ILFD). It is the natural extension of the analysis 
performed in~\cite{BDG3}, where we analytically proved the experimental
and numerical results contained in~\cite{OBK,OBYK} by providing
explicit formulae for the width of the plateaux appearing in
the devil's staircase. More precisely in~\cite{BDG3}
we found the following result.
Denote by $\omega$ and $\Omega$ the frequencies of the driving signal
and of the output signal of the ILFD, respectively, with
$\mu$ the driving amplitude. Then, if for $\rho\in\QQQ$ we call
$\Delta\omega(\rho)=\{\omega:\omega/\Omega=\rho\}$ the width of the
corresponding locking interval, we showed that
$\Delta\omega(\rho)$ satisfies $\Delta\omega(2n/k) =O(\mu^{k})$
and $\Delta\omega((2n-1)/k) = O(\mu^{2k})$ for all $k,n \in \NNN$
such that $2n/k$ and $(2n-1)/k$ are irreducible fractions. In particular
this implies that the largest plateaux correspond to even integer
values of the ratio $\omega/\Omega$.

In this paper we have extended the above results: we studied
the system of ordinary differential equations (\ref{eq:2.1}),
(\ref{eq:2.2}), (\ref{eq:2.3}), which describe the ILFD, with a
more general driving term in the form of any analytic
periodic function (we confined ourselves to functions containing 
odd-harmonics only in order to make the analysis more transparent
and yet without any significant loss of generality). In~\cite{BDG3}
we used $f(t)= \sin(t)$ (one harmonic only), as in~\cite{OBK,OBYK}.
Here, we studied the locking intervals $\Delta\omega(\rho)$
by using in~(\ref{eq:2.3}) a $2\pi$-periodic function
of the form $f(t) = \sum_{\nu=1}^{\io} \hat f_{\nu} \sin \nu t$, with
$|\hat f_{\nu} | \le \Phi\,{\rm e}^{-\xi |\nu|}$ (by analyticity).
We found that, for any $\rho=p/q\in\QQQ$, with $p,q$
relatively prime integers, the key condition for the existence of
the locking $\omega=\rho\Omega$ (and hence of a plateau),
is that there exists $\nu$ such that $\hat f_{\nu}\neq0$ and
$2|\nu'|q=|\nu|p$  and some $\nu'\in\ZZZ$. This condition is
certainly satisfied if, for instance, $|\nu'|=p$ and $|\nu|=2q$,
provided $\hat f_{2q} \neq 0$, or $|\nu'|=2p$ and $|\nu|=q$,
provided $\hat f_{q} \neq 0$. Thus, for any resonance $p\!:\!q$
one has a plateau $\Delta\omega(p/q)$ which to first order is given
by (\ref{eq:4.1}) and (\ref{eq:4.6}). In particular,
to leading order, the width of the \Arnold\ tongues is expressed as
$\Delta\omega_{1}(\rho) \approx 2\mu\rho\Omega_{0}|\bar r_{1}|^{-1}
|\hat f_{\nu_{0}} | \, | \overline K_{\nu_{0}\rho}(\rho) |,$
where $\rho=p/q$, $\Omega_{0}$ and $\bar r_{1}$
are constants depending on the unperturbed system (but not on
the driving) and $\nu_{0} \geq 1$ denotes the integer which provides 
the leading coefficient in the sum (\ref{eq:4.6}).
Note that the formula reduces to the one obtained in~\cite{BDG3} 
--- as it should --- if $f(t)= \sin(t)$: in that case $\hat f_{\nu}
\neq 0$ only for $\nu=1$, so that $q=1$ and $p\in 2\NNN$.

Moreover, by keeping in (\ref{eq:4.6}) the whole sum, we obtain
$\left| \Delta\omega(\rho) \right| \le \mu
C\, p^{2}q^{-1} {\rm e}^{-\xi_{1} p} {\rm e}^{-\xi q}$
where $C$ is a constant independent of $p$ and $q$, thereby
showing that the \Arnold\ tongues have width proportional to $\mu$,
but with  proportionality constants which decay exponentially
with $p$ and $q$.

We have also computed analytically the contribution of the second order,
namely the coefficient of $\mu^2$. In this case one needs
to compute the first order solution $(u_1(\tau), \dot u_1(\tau))$,
with $u_1(\tau)$ given in (\ref{eq:4.17}), which rather complicates
the analysis. We found $\Delta\omega(\rho) =
\mu^{2} \Delta_{2}\omega(\rho) + O(\mu^{3}),$
which replaces~(\ref{eq:4.1}) when the first order vanishes.
For instance, if $f$ contains only the first harmonics then
the condition for locking onto a $p\!:\!q$ resonance
becomes: $2|\nu'|q=|\nu_{1}+\nu_{2}|p$, with $\hat f_{\nu_{1}}
\hat f_{\nu_{2}} \neq 0,$ is satisfied for some $\nu'\in\ZZZ$.
This shows that when $f(t)=\sin t$, as in~\cite{OBK,OBYK,BDG3},
the plateaux corresponding to odd $\rho$ are of order $\mu^{2}$. 

Higher order contributions can in principle be computed with a very
similar strategy (see (\ref{eq:4.23}) and (\ref{eq:4.24}));
the important point to notice is that to all order $k$
the coefficients $\Delta_{k}\omega(\rho)$ decay exponentially
in both $p$ and $q$. Naturally higher order terms become
dominant when all the terms of smaller order vanish.

To complete our investigation, we computed the functions
$\gotD_{1}(\tau_{0})$ and $\gotD_{2}(\tau_{0})$ numerically,
from which the tongue widths $\Delta\omega(\rho)$ and
$\Delta\omega_2(\rho)$ can be calculated, via equations~(\ref{eq:4.1})
and~(\ref{eq:4.22}) respectively. Some of the techniques required
to carry out this computation are described in Section~\ref{sec:4}.
We then computed a set of \Arnold\ tongues, which was sufficiently
large for testing the numerics on the basis of the theoretical
predictions. In particular, we computed the width of the tongues for
two types of forcing: (i) only one harmonic and (ii) all harmonics
present in the Fourier expansion. Our computational results are 
in excellent alignment with the theory as stated above,
which supports our belief that the locking charts in
figure~\ref{fig:tongues} are accurate. These two charts clearly
demonstrate the dominance of the 2:1 and the 4:1 resonances.
Furthermore, a comparison indicates that the location of the tongues
is robust under generic perturbations; the differences in the shapes
of the tongues are small.

\appendix

\zerarcounters
\section{Error in the interpolation scheme}
\label{app:a}

Starting from equation~(\ref{I_of_t}), we expand the sine
functions in terms of complex exponentials to obtain
\begin{equation}
I_K(t) = \frac{1}{K}\sum_{|j|\leq (K-1)/2} {\rm e}^{2 {\rm i} j\pi t}.
\label{ikc}
\end{equation}
Setting $t = \tau/T_0$, so that the scheme can be used to interpolate
a periodic function $x(t)$ of arbitrary period $T_0$ in terms of $\tau$,
we have the Fourier expansion
$$ x(t) = \sum_{n\in\ZZZ}\alpha_n {\rm e}^{2 {\rm i} n\pi t}. $$
We interpolate $x(t)$ by
\begin{equation}
\hat{x}(t) = \sum_{j = 0}^{K-1} x(j/K) I_K(t - j/K) = \sum_{|m| \leq
(K-1)/2}\beta_m {\rm e}^{2 {\rm i} m\pi t} ,
\label{hatx}
\end{equation}
where the last inequality follows from equation~(\ref{ikc}). In order to
determine how well $x(t)$ is approximated by $\hat{x}(t)$, we need to
compare $\alpha_m$ with $\beta_m$. Substituting for $x(j/K)$
and $I_K(t)$ in equation~(\ref{hatx}) and rearranging, we find
$$ \hat{x}(t) = \sum_{|m|\leq (K-1)/2} {\rm e}^{2 {\rm i} m\pi
t}\left\{\sum_{n\in\ZZZ}\frac{\alpha_n}{K}\sum_{j=0}^{K-1}
{\rm e}^{2 {\rm i} (n-m)j\pi/K}\right\} , $$
where the term in braces is equal to $\beta_m$. The sum over $j$
is equal to $[1 - {\rm e}^{2{\rm i}(n-m)\pi}]/[1 -
{\rm e}^{2{\rm i}(n-m)\pi/K}]$
provided that $(n-m)\neq pK$, $p\in\ZZZ$,
and is equal to $K$ otherwise. Hence,
$$ \frac{1}{K}\sum_{j=0}^{K-1} {\rm e}^{2{\rm i}(n-m)j\pi/K} = \left\{
\begin{array}{l@{\hskip 0.2in}l}
0 & n-m\neq p K\\
1 & n-m = p K.
\end{array}\right. $$
Hence
$$ \beta_m = \sum_{p\in\ZZZ}\alpha_{m + pK} = \alpha_m + \alpha_{m-K} +
\alpha_{m+K} + \ldots $$
Now, since $x(t)$ is the solution of an ODE with analytic coefficients,
it is itself analytic, and so, for all $n$,
$|\alpha_n| < C_1 {\rm e}^{-C_2 |n|}$, where $C_1, C_2$
are positive real constants. Thus, 
\begin{equation}
|\beta_m - \alpha_m| < C_1 {\rm e} ^{-C_2 K}
\nonumber \\ 
\end{equation}
and hence, by choosing $K$ sufficiently large, the interpolation error
can be made as small as we please.

\vspace{.5truecm}
\noindent \textbf{Acknowledgments.} We thank Peter Kennedy for 
bringing this problem to our attention, and for the interest shown
in our results. We are also indebted to
Bill Christmas for useful discussions on interpolation.
FS was supported by the Engineering and Physical Sciences Research
Council (EPSRC) Grant no.\ EP/D063906/1.



\begin{thebibliography}{99}

{\small

\bibitem{AC}{
A.A. Abidi, L.O. Chua,
\textit{On the dynamics of Josephson-junction circuits},
Electron. Circuits Syst.
\textbf{3} (1979), no. 4, 974--980. }
%
\bibitem{AMOQR}{
A. Amann, M.P. Mortell, E.P. O'Reilly, M. Quinlan, D. Rachinskii,
\textit{Mechanisms of syncronization in frequency dividers},
IEEE Trans. Circuits Syst. I. Regul. Pap.
\textbf{56} (2009), no.1, 190--199. }
%
\bibitem{A}{
V.I. \Arnold ,
\textit{Geometrical methods in the theory of
ordinary differential equations},
Grund\-lehren der Mathematischen Wissenschaften Vol. 250,
Springer, New York, 1988. }
%
\bibitem{BDG3}{
M. Bartuccelli, J. Deane, G. Gentile,
\textit{Frequency locking in the injection-locked frequency divider equation},
Proc. R. Soc. Lond. Ser. A Math. Phys. Eng. Sci.
\textbf{465} (2009), no. 2101, 283--306. }
%
\bibitem{Cha}{
Y. F. Chang and G. Corliss, 
\textit{Ratio-like and recurrence relation tests
for convergence of series},
J. Inst. Math. Appl.
\textbf{25} (1980), 349--359. }
%
\bibitem{C}{
W.A. Coppel,
\textit{Some quadratic systems with at most one limit cycle},
Dynamics reported Vol. 2, 61--88,
Dynam. Report. Ser. Dynam. Systems Appl., 2, Wiley, Chichester, 1989. }
%
\bibitem{DDK}{
S. Daneshgar, O. De Feo, M.P. Kennedy,
\textit{Obsrvations concerning the locking range in a complementary
differential $LC$ injection-locked frequency divider.
Part I: qualitative analysis},
IEEE Trans. Circuits Syst. I. Regul. Pap., to appear. }
%
\bibitem{GBD3}{
G. Gentile, M. Bartuccelli, J. Deane,
\textit{Bifurcation curves of subharmonic solutions and
Melnikov theory under degeneracies},
Rev. Math. Phys.
\textbf{19} (2007), no. 3, 307--348. }
%
\bibitem{GH}{
J. Guckenheimer, Ph. Holmes,
\textit{Nonlinear oscillations, dynamical systems,
and bifurcations of vector fields},
Applied Mathematical Sciences Vol. 42,
Springer, New York, 1990. }
%
\bibitem{H}{
Ph. Hartman,
\textit{Ordinary differential equations},
Classics in Applied Mathematics Vol. 38,
Society for Industrial and Applied Mathematics (SIAM),
Philadelphia, PA, 2002. }
%
\bibitem{He}{
M.R. Herman,
\textit{Mesure de Lebesgue et nombre de rotation},
Geometry and topology (Proc. III Latin Amer. School of Math.,
Inst. Mat. Pura Aplicada CNPq, Rio de Janeiro, 1976),
pp. 271--293, Lecture Notes in Mathematics, Vol 597,
Springer, Berlin, 1977. }
%
\bibitem{KC}{
M.P. Kennedy, L.O. Chua,
\textit{Van der Pol and chaos},
IEEE Trans. Circuits Syst. I. Regul. Pap.
{\bf 33} (1986), no. 10, 974--980. }
%
\bibitem{L1}{
M. Levi,
\textit{Nonchaotic behavior in the Josephson junction},
Phys. Rev. A (3)
\textbf{37} (1988), no. 3, 927--931. } 
%
\bibitem{M}{
P. Maffezzoni,
\textit{Analysis of oscillator injection locking
through phase-domain impulse-response},
IEEE Trans. Circuits Syst. I. Regul. Pap.
\textbf{55} (2008), no. 5, 1297--1305. }
%
\bibitem{OBK}{
D. O'Neill, D. Bourke, M.P. Kennedy,
\textit{The Devil's staircase as a method of comparing injection-locked
frequency divider topologies},
Proceedings of the 2005 European Conference on
Circuit Theory and Design, 2005, Vol. III, pp. 317--320. }
%
\bibitem{OBYK}{
D. O'Neill, D. Bourke, Zh. Ye, M.P. Kennedy,
\textit{Accurate modeling and experimental validation of
an injection-locked frequency divider},
Proceedings of the 2005 European Conference on
Circuit Theory and Design, 2005, Vol. III, pp. 409--412. }
%
\bibitem{PL}{
U. Parlitz, W. Lauterborn,
\textit{Period-doubling cascades and devil's staircases
of the driven van der Pol oscillator},
Phys. Rev. A
\textbf{36} (1987), no. 3, 1428 - 1434. }
%
\bibitem{PZC}{
L. Pivka, A.L. Zheleznyak, L.O. Chua,
\textit{\Arnold\ tongues, devil's staircase, and self-similarity
in the driven Chua's circuit},
Internat. J. Bifur. Chaos Appl. Sci. Engrg.
\textbf{4} (1994), no. 6, 1743--1753. }
%
\bibitem{numrec}{
W.H. Press, S.A. Teukolsky, W.T. Vetterling, B.P. Flannery, 
\textit{Numerical Recipes in C} (2nd edition, 1992),
Cambridge University Press, Cambridge. }
%
\bibitem{QWZ}{
M. Qian, J.-Z. Wang, X.-J. Zhang,
\textit{Resonant regions of Josephson junction equation in case
of large damping},
Phys. Lett. A
\textbf{372} (2008), no. 20, 3640-3644. }
%
\bibitem{FS}{
F. Schilder, B. B. Peckham,
\textit{Computing \Arnold\ tongue scenarios},
J. Comput. Phys.
\textbf{220} (2006) no. 2, 932--951. }
%
\bibitem{PM}{
B. van der Pol, J. van der Mark,
\textit{Frequency demultiplication},
Nature
\textbf{120} (1927), no. 3019, 363--364. }
%
\bibitem{YXK}{
Zh. Ye, T. Xu, M.P. Kennedy,
\textit{Locking range analysis for injection-locked frequency dividers},
Proceedings of the 2006 IEEE International Symposium on
Circuits and Systems, 2006, 4070--4073. }
%
\bibitem{Z}{
Zh. F. Zhang,
\textit{Proof of the uniqueness theorem of limit cycles
of generalized Li\'enard equations},
Appl. Anal.
\textbf{23} (1986), no. 1-2, 63--76. }

}

\end{thebibliography}
\end{document}